\documentclass[11pt,a4paper]{article}
\usepackage[T1]{fontenc}
\usepackage{amsmath,amssymb,amsthm}
\usepackage{natbib}
\usepackage[margin=1in,a4paper]{geometry}
\usepackage{booktabs, multicol, multirow}
\usepackage[colorlinks,citecolor=blue,urlcolor=blue]{hyperref}
\usepackage{mathtools}
\mathtoolsset{showonlyrefs=true}
\usepackage{subcaption}
\usepackage{float}
\usepackage[linesnumbered,ruled,vlined]{algorithm2e}
\usepackage{tikz}
\usetikzlibrary{shapes,arrows}
\usepackage{caption}

\numberwithin{equation}{section}
\newtheorem{proposition}{Proposition}

\newtheorem{assumption}{Assumption}
\newtheorem{theorem}{Theorem}
\theoremstyle{remark}
\newtheorem{remark}{Remark}
\theoremstyle{remark}

\allowdisplaybreaks[4]

\title{Simulation of Multidimensional Diffusions with Sticky Boundaries via Markov Chain Approximation}
\author{Christian Meier\thanks{Department of Systems Engineering and Engineering Management, The Chinese University of Hong Kong, Hong Kong. Email: meier@se.cuhk.edu.hk.}\and Lingfei Li\thanks{Corresponding author. Department of Systems Engineering and Engineering Management, The Chinese University of Hong Kong, Hong Kong. Email: lfli@se.cuhk.edu.hk.} \and Gongqiu Zhang\thanks{School of Science and Engineering, The Chinese University of Hong Kong, Shenzhen, China. Email: zhanggongqiu@cuhk.edu.cn.}}

\begin{document}

\maketitle
\begin{abstract}
	We develop a new simulation method for multidimensional diffusions with sticky boundaries. The challenge comes from simulating the sticky boundary behavior, for which standard methods like the Euler scheme fail. We approximate the sticky diffusion process by a multidimensional continuous time Markov chain (CTMC), for which we can simulate easily. We develop two ways of constructing the CTMC: approximating the infinitesimal generator of the sticky diffusion by finite difference using standard coordinate directions, and matching the local moments using the drift and the eigenvectors of the covariance matrix as transition directions. The first approach does not always guarantee a valid Markov chain whereas the second one can. We show that both construction methods yield a first order simulation scheme, which can capture the sticky behavior and it is free from the curse of dimensionality. We apply our method to two applications: a multidimensional Brownian motion with all dimensions sticky which arises as the limit of a queuing system with exceptional service policy, and a multi-factor short rate model for low interest rate environment in which the stochastic factors are unbounded but the short rate is sticky at zero. 
	
	
	\hspace{1.95cm}Monte Carlo simulation.
	
\end{abstract}


\section{Introduction}
\label{sec:introduction}
Diffusion processes with sticky boundaries arise as natural models in various fields like biology, physics, queuing systems and finance.
Sticky boundary behavior is discovered by \cite{feller1952} in the one-dimensional case and a historical account is given by \cite{peskir2015}.
Generalizations to multiple dimensions are developed in \cite{wentzell1959,wentzell1960}.
The existence and uniqueness of multidimensional stochastic differential equations with sticky boundaries is studied in \cite{graham1988} and \cite{ikeda1989}. In addition,  \cite{grothaus2017} and \cite{racz2015} provide detailed analysis of multidimensional sticky Brownian motions. Some applications of sticky diffusions can be found in e.g., \cite{graham1989}, \cite{kalda2007}, \cite{racz2015}, \cite{fattler2016} and \cite{nie2020}. 

The present paper focuses on the computational aspect of multidimensional diffusions with sticky boundaries, which is not well studied in the literature. Under certain conditions, the value function for the expectation of the process is the solution to a parabolic PDE with Wentzell boundary condition. One can try to solve this type of PDE numerically by finite difference or finite element method; see e.g., \cite{kolkovska2007}, \cite{kovacs2017}, \cite{bansch2020} and \cite{gander2021}. Nevertheless, the numerical PDE approach suffers from the curse of dimensionality, making it computationally intractable for high-dimensional problems. Furthermore, simulating sample paths from the model is important in applications, which numerical PDE methods cannot do. 

Recently, \cite{meier2021} propose a computational method for general one-dimensional diffusions with sticky boundaries using continuous time Markov chain (CTMC) approximation (also see \cite{bou-rabee2020} for one-dimensional sticky Brownian motions). Under their scheme, the value function can be computed efficiently and accurately, and they simulate paths from the CTMC to approximate the original model, which can capture the sticky behavior. The present paper extends this approach to multidimensional diffusions with sticky boundaries in nontrivial ways and our goal is to develop a simulation method so that we can generate paths of the process and calculate the value function using the Monte Carlo method. By virtue of the Feynman-Kac Theorem, one can also apply our simulation method to solve high-dimensional parabolic PDEs with Wentzell boundary condition. 

Markov chain approximation has been developed in the literature for some multidimensional diffusions. See \cite{kushner2001} for a comprehensive study on constructing discrete time Markov chains for diffusions without boundaries or with reflecting boundaries to solve stochastic control problems. \cite{kirkby2020} consider pricing multi-asset financial options under the multidimensional geometric Brownian motion (GBM) model. They apply a suitable change of variables to decorrelate the multidimensional Brownian motion and construct a CTMC to approximate each coordinate, which is independent of each other after the transform. However, the success of their approach depends on the specific forms of the drift vector and diffusion matrix of the GBM, which cannot be applied to general diffusions. To price financial options in coupled two-dimensional diffusions models, \cite{xi2019simultaneous} construct CTMC approximation in a way that is equivalent to approximating the generator by finite difference. This approach has some drawbacks as we will discuss below. Moreover, they don't consider diffusions with boundary behaviors. \cite{cui2018} develop a two-layer CTMC approximation for pricing options in stochastic local volatility models which form a class of two-dimensional diffusions. Under their approximation, the stock price is approximated by a regime-switching CTMC. \cite{cui2020efficient} further considers how to simulate these processes. They simulate the one-dimensional variance process by a CTMC that approximates it and then sample the integrated variance conditioned on the start and end points of the variance process using a Fourier sampler for the CTMC variance model. Finally, there are also various papers on CTMC approximation of one-dimensional Markov processes with applications in finance; see e.g., \cite{mijatovic2013}, \cite{cai2015},  \cite{eriksson2015}, \cite{cui2018single}, \cite{li2018error}, \cite{zhang2019,zhang2021parisian,zhang2021drawdown,zhang2021nonsmooth}and \cite{zhang2021pricing}.

In this paper, we develop CTMC approximation for general multidimensional diffusions with sticky boundaries. Although our focus is on sticky diffusions, our approach can be adapted to construct CTMC approximation for diffusions with other types of boundary behaviors. Below we discuss two major issues which would make the construction nontrivial.

\smallskip
\noindent (1) The first issue is what directions should be used to move the CTMC. In the one-dimensional case, the CTMC can only move to the left or right along the real line. However, in a multidimensional setting, there are many possible directions and it is not clear at all what directions should be used. A standard idea is to approximate the differential operators in the infinitesimal generator of the sticky diffusion by finite difference, which approximates partial derivatives of the kind $\partial_{x^ix^j}$ using changes along the $i$th and $j$th dimensions. This approximation implies that the CTMC can only have at most two of its coordinates moving in every transition, which is not realistic if multiple coordinates are strongly correlated in the original model. Furthermore, as we will show later, the transition rates obtained from the finite difference approach can be negative, and hence they are not valid transition rates for a CTMC. 

In this paper, we propose to use the drift vector and the eigenvectors of the covariance matrix, which are orthogonal to each other, as the directions to move the CTMC. We will simply call it as the eigendecomposition approach. Using eigendirections allows multiple coordinates to move together in one transition, which can capture the strong correlations among components of the diffusion part. We obtain the transition rates by matching the local moments of the drift and the diffusion parts. Unlike the finite difference approach, the resulting transition rates of the eigendecomposition approach are always valid for a CTMC. We will also show that if the same step size is used, the eigendecomposition approach can significantly reduce the CTMC approximation error compared with the finite difference approach, and this is because the set of eigendirections contains those ones along which the process varies the most. 

\smallskip
\noindent (2) The second issue is how to deal with the sticky boundary. We adjust the step size of the CTMC carefully to avoid moving out of the boundary when the process is close to it. Moreover, we must correctly capture the sticky behavior on the boundary, and this is achieved by matching the local moments of the CTMC with the diffusion on the boundaries. 
 
Euler scheme is the standard method for simulating multidimensional diffusions without boundaries. If a diffusion has an absorbing or reflecting boundary, modifications can be made in the Euler scheme to correctly simulate these behaviors (see e.g., \cite{gobet2000}, \cite{gobet2001}, \cite{bossy2004}, \cite{bayer2010}, \cite{gobet2010}, \cite{nystrom2010}). However, it is unclear how to modify the Euler scheme to simulate the sticky boundary behavior and its failure in the one-dimensional case is documented in \cite{meier2021}. Our approach provides a convergent simulation scheme that fills the gap in the literature for the sticky boundary behavior. We will prove that the weak convergence order of our method is one in terms of the step size of the CTMC. This is comparable to the convergence rate of the Euler scheme which is first order in the time step. 

The rest of this paper is organized as follows. Section \ref{sec:characterization_multidimensional_sticky_processes} characterizes multidimensional diffusion processes with sticky boundaries. In particular, we derive an alternative SDE for these processes which is easier to simulate and their infinitesimal generator. In Section \ref{sec:simulation_2d_2_sticky_dimensions} we show how to construct a multidimensional CTMC using the finite difference approach and the eigendecomposition approach. We further prove the weak convergence order of our scheme under both approaches is 1. Section \ref{sec:numerical_examples} demonstrates the performance of these two construction approaches on two applications and also compares discrete and exact simulation of CTMCs. Finally, Section \ref{sec:conclusion} concludes. The appendix provides proofs for all the results. 

To close this section, we introduce some notations. Throughout the paper, for a vector $x=(x^1,x^2,\ldots,x^d)$, we use superscripts for the coordinates and we write the power of a coordinate as e.g.,  $\left(x^1\right)^2$. In addition, we use $\Vert\cdot\Vert$ for the Euclidean norm of a vector.

\section{Multidimensional Diffusions with Sticky Boundaries}
\label{sec:characterization_multidimensional_sticky_processes}

Consider a general $d$-dimensional diffusion. Some or all of its dimensions have a sticky boundary, and the others (if any) are unbounded. The diffusion lives on $\bar{\mathbb{S}}=\mathbb{S}\cup\partial\mathbb{S}$, where
\begin{equation}
	\mathbb{S}=\{x\in\mathbb{R}^d:\Phi\left(x\right)>0\}, \qquad\partial\mathbb{S}=\{x\in\mathbb{R}^d:\Phi\left(x\right)=0\},
\end{equation}
for some $\Phi\in C_b^2(\mathbb{R}^d)$, which is the space of twice continuously differentiable functions where the function itself, first and second order derivatives are bounded. Let $n(x)$ be the unit normal vector at the boundary point $x$ pointing inwards, which is given by 
$n=\nabla\Phi(x)/{\Vert\nabla\Phi(x)\Vert}$.

\begin{remark}
	In the applications in Section \ref{sec:numerical_examples}, we have $\mathbb{S}=\{x\in\mathbb{R}^d:x^1,\ldots,x^{\hat{d}}>0\}$ where $\hat{d}$ is the number of dimensions exhibiting stickiness.
	A choice of $\Phi$ in this case is given by (see Remark 3 in Section 1, \cite{graham1988})
	\begin{equation}\label{eq:definition_Phi_orthant}
		\Phi\left(x\right)=\prod_{i=1}^{\hat{d}}\left(1-\exp\left(-x^i\right)\right).
	\end{equation}
The boundary for the $i$-th coordinate is given by $\{x\in\mathbb{R}^d: x^i=0\}$ for $i=1,\cdots,\hat{d}$.
\end{remark}

\subsection{The Formulation}
We specify the sticky diffusion as the solution to some stochastic differential equation (SDE). Consider the following measurable functions:
\begin{align}
	\mu&:\bar{\mathbb{S}}\to\mathbb{R}^d,\,\,\,\,\,\Sigma:\bar{\mathbb{S}}\to\mathbb{R}^{d\times d}, \\
	\beta&:\partial\mathbb{S}\to\mathbb{R}^{d},\ \Gamma:\partial\mathbb{S}\to\mathbb{R}^{d\times d},\  \rho:\partial\mathbb{S}\to\mathbb{R}_{+},
\end{align}
Let $A=\Sigma\Sigma^\top$ and $G=\Gamma\Gamma^\top$. We consider the following SDE:
\begin{equation}\label{eq:multidimensional_sticky_SDE}
	dX_t=\mu\left(X_t\right)\left(dt-\rho\left(X_t\right)dL_t\right)+\Sigma\left(X_t\right)dM_t+\beta\left(X_t\right)dL_t+\Gamma\left(X_t\right)dB_{L_t}
\end{equation}
with $X_0=x$. Here, $B$ is a $d$-dimensional standard Brownian motion, $L$ is the local time process of $X$ on the boundary, i.e., 
\begin{equation*}
	L_t=\int_0^tI\left(X_s\in\partial\mathbb{S}\right)dL_s,
\end{equation*}
and $M$ is a $d$-dimensional continuous martingale with
\begin{equation*}
	\langle M^i,M^j\rangle_t=\delta_{i,j}\left(t-\int_0^t\rho\left(X_s\right)dL_s\right),\ i,j=1,\ldots,d
\end{equation*}
where $\delta_{i,j}$ is the Kronecker delta. We make the following assumption. 
\begin{assumption}\label{assumptions:existence}
	Assume the following properties hold. 
	\begin{enumerate}
		\item $\mu(x)$ and $\Sigma(x)$ are bounded and Lipschitz continuous and $\Sigma(x)$ has full rank on $\bar{\mathbb{S}}$, and $n(x)^\top A(x) n(x)>C>0$ on $\partial{\mathbb{S}}$ for some constant $C$.
		\item $\beta(x)$ and $\Gamma(x)$ are bounded and Lipschitz continuous on $\partial\mathbb{S}$.
		\item $\Gamma(x)^\top n(x)=0$.
		\item $\rho(x)$ is strictly positive and there exists some constant $C>0$ such that $\beta(x)^\top n(x)>C$ for any $x\in\partial\mathbb{S}$.
	\end{enumerate}
\end{assumption}

Under Assumption \ref{assumptions:existence}, we can apply Theorems I.13 in \cite{graham1988} to conclude that there exists a unique weak solution to \eqref{eq:multidimensional_sticky_SDE}. Furthermore, it follows that 
the sojourn condition holds, which is
\begin{equation}\label{eq:sojourn_condition_multidimensional}
	I\left(X_t\in\partial\mathbb{S}\right)dt=\rho\left(X_t\right)dL_t.
\end{equation}

\begin{remark}
	We explain the implication of Condition 3 in Assumption \ref{assumptions:existence}. Suppose $\Phi$ is given by \eqref{eq:definition_Phi_orthant} and $\mathbb{S}=\{x\in\mathbb{R}^d:x^1,\ldots,x^{\hat{d}}>0\}$. Consider a point $x$ such that $x^i=0$ and $x^j\neq 0$ for $j\neq i$. Then $n(x)=e_i$ in which the $i$th coordinate is one and all others are zero.  The condition $\Gamma(x)^\top n(x)=0$ implies that the $i$th row of $\Gamma(x)$ is zero. Hence the $i$th row and $i$th column of $G(x)=\Gamma(x)\Gamma(x)^\top$ are zero. In other words, there is no diffusion for the $i$th coordinate.  At the point $x$ where $x^1=\cdots=x^{\hat{d}}=0$, similarly one obtains that the first to $\hat{d}$th rows and columns of $G(x)$ are zero. So there is no diffusion for the first $\hat{d}$ coordinates at this point. In general, there is no diffusion for a coordinate that shows sticky boundary behavior when it is at zero. 
\end{remark}

\begin{remark}
Conditions 1 and 2 on the boundedness and Lipschitz continuity of these functions are standard conditions in the literature to ensure well-posedness for \eqref{eq:multidimensional_sticky_SDE}, but they are not necessary for our simulation method to work. For a specific SDE, one may be able to show existence and uniqueness of a weak solution under weaker conditions. Condition 4 guarantees that the drift vector is pointing in the right direction to avoid taking the process out of the state space. 
\end{remark}

Simulating the SDE \eqref{eq:multidimensional_sticky_SDE} directly requires the simulation of the local time process. We refer readers to some references on this topic; see \cite{etore2013,etore2018} and \cite{bourza2020} for one-dimensional processes and \cite{blanchet2018} for multidimensional reflected Brownian motion. However, we cannot apply their methods in our problem which involves a Brownian motion time changed by a local time process and they are not independent.  In our approach, we will rewrite \eqref{eq:multidimensional_sticky_SDE} to get rid of the local time term by using the sojourn condition \eqref{eq:sojourn_condition_multidimensional}.
Firstly, one can see that
\begin{equation*}
	dt-\rho\left(X_t\right)dL_t=dt-I\left(X_t\in\partial\mathbb{S}\right)dt=\left(I\left(X_t\in\bar{\mathbb{S}}\right)-I\left(X_t\in\partial\mathbb{S}\right)\right)dt=I\left(X_t\in\mathbb{S}\right)dt.
\end{equation*}
Secondly, the quadratic variation process of $M$ can be rewritten as 
\begin{equation*}
	\langle M^i,M^j\rangle_t=\delta_{i,j}\int_0^t\left(I\left(X_s\in\bar{\mathbb{S}}\right)-I\left(X_s\in\partial\mathbb{S}\right)\right)ds=\delta_{i,j}\int_0^t I\left(X_s\in\mathbb{S}\right)ds.
\end{equation*}
Application of the martingale representation theorem (see e.g., Theorem 4.2 in \cite{karatzas1991}) shows that $M$ can be represented as (after possibly enlarging the underlying probability space)
\begin{equation*}
	M_t=\int_0^tI\left(X_s\in\mathbb{S}\right)dB_{1,s},
\end{equation*}
where $B_{1,s}$ is a standard $d$-dimensional Brownian motion independent of $B_t$ (see Remark 1 in Section 3 of \cite{graham1988}).
Thirdly, we have 
\begin{equation*}
	\beta\left(X_t\right)dL_t=\beta\left(X_t\right)\frac{1}{\rho\left(X_t\right)}I\left(X_t\in\partial\mathbb{S}\right)dt,
\end{equation*}
which is well-defined by the positivity of $\rho$.
Lastly, another application of the martingale representation theorem yields
\begin{equation*}
	B_{L_t}=\int_0^t\sqrt{\frac{dL_s}{ds}}dB_{2,s}=\int_0^t\frac{1}{\sqrt{\rho\left(X_s\right)}}I\left(X_s\in\partial\mathbb{S}\right)dB_{2,s},
\end{equation*}
where $B_{2,t}$ is another standard $d$-dimensional Brownian motion independent of $B_{1,t}$. Putting them together, we obtain
\begin{align}
	dX_t&=\mu\left(X_t\right)I\left(X_t\in\mathbb{S}\right)dt+\Sigma\left(X_t\right)I\left(X_t\in\mathbb{S}\right)dB_{1,t} \nonumber \\
	&\qquad+\hat{\beta}\left(X_t\right)I\left(X_t\in\partial\mathbb{S}\right)dt+\hat{\Gamma}\left(X_t\right)I\left(X_t\in\partial\mathbb{S}\right)dB_{2,t},\label{eq:multidimensional_sticky_SDE_adjusted}
\end{align}
with
\begin{equation}\label{eq:definition_hat_beta}
	\hat{\beta}\left(x\right)=\frac{\beta\left(x\right)}{\rho\left(x\right)},\qquad\hat{\Gamma}\left(x\right)=\frac{\Gamma\left(x\right)}{\sqrt{\rho\left(x\right)}},\qquad\hat{G}\left(x\right)=\hat{\Gamma}\left(x\right)\hat{\Gamma}\left(x\right)^\top\qquad\text{for}\ x\in\partial\mathbb{S}.
\end{equation}
One can see from this alternative expression of the SDE that $\mu$ and $\Sigma$ describe the evolution of $X$ in the interior of the state space, whereas $\hat{\beta}$ and $\hat{\Gamma}$ describe the behavior on the boundary. Hereafter, we will work with the formulation 
\eqref{eq:multidimensional_sticky_SDE_adjusted} for the SDE.

\subsection{The Transition Semigroup and the Infinitesimal Generator}
Under Assumption \ref{assumptions:existence}, $X$ is a Feller process (Theorem 5.11, \cite{dynkin1965}). Consider the transition operator defined as
\begin{equation}\label{eq:feynman_kac_operator_def}
	\mathcal{P}_tf\left(x\right)=\mathbb{E}_x\left(f\left(X_t\right)\right),\qquad x\in\bar{\mathbb{S}}.
\end{equation}
Then $(\mathcal{P}_t)_{t\ge0}$ is a strongly continuous semigroup of contractions on $C_0(\bar{\mathbb{S}})$, the space of continuous functions on $\bar{\mathbb{S}}$ vanishing at $\infty$. The infinitesimal generator of the process is defined as 
\begin{equation}\label{eq:generator-def}
	\mathcal{G}f=\text{s-}\underset{t\searrow 0}{\lim}\ \frac{\mathcal{P}_tf-f}{t},
\end{equation}
for functions where the limit exists and $\text{s-}$ indicates the limit is taken under the norm of $C_0(\bar{\mathbb{S}})$. We next show the formula for the generator. Introduce two operators:
\begin{align}
	\mathcal{A}f\left(x\right)&=\sum_{i=1}^{d}\mu^i\left(x\right)\frac{\partial}{\partial x^i}f\left(x\right)+\frac{1}{2}\sum_{i,j=1}^{d}A^{i,j}\left(x\right)\frac{\partial^2}{\partial x^i\partial x^j}f\left(x\right)\\
	&=\left(\partial_{x}f\right)^\top\mu\left(x\right)+\frac{1}{2}\textrm{Tr}\left(\Sigma\left(x\right)^\top\left(\partial_{xx}f\right)\Sigma\left(x\right)\right),
	\qquad x\in\bar{\mathbb{S}}, \label{eq:definition_A_operator} 
\end{align}
and
\begin{align}
	\mathcal{K}f\left(x\right)&=\sum_{i=1}^{d}\hat{\beta}^{i}\left(x\right)\frac{\partial}{\partial x^{i}}f\left(x\right)+\frac{1}{2}\sum_{i,j=1}^{d}\hat{G}^{i,j}\left(x\right)\frac{\partial^2}{\partial x^i\partial x^j}f\left(x\right)\\
	&=\left(\partial_{x}f\right)^\top\hat{\beta}\left(x\right)+\frac{1}{2}\textrm{Tr}\left(\hat{\Gamma}\left(x\right)^\top\left(\partial_{xx}f\right)\hat{\Gamma}\left(x\right)\right),\qquad x\in\partial\mathbb{S}, \label{eq:definition_K_operator}
\end{align}
where $\partial_xf$ is the gradient of $f$, $\partial_{xx}f$ is its Hessian matrix and $\textrm{Tr}(\cdot)$ denotes the trace operator.

\begin{theorem}\label{th:generator_multidimensional_sticky_diffusion}
	Under Assumption \ref{assumptions:existence}, for any $f\in C^2_0(\bar{\mathbb{S}})$ (functions in $C_0(\bar{\mathbb{S}})$ with their first and second order derivatives continuous and vanishing at infinity) that satisfies
	\begin{equation}\label{eq:Wentzell}
		\mathcal{A}f(x)=\mathcal{K}f(x),\ \text{for any}\ x\in \partial\mathbb{S},
	\end{equation} 
	the limit \eqref{eq:generator-def} exists and
	\begin{equation}\label{eq:generator}
		\mathcal{G}f\left(x\right)=I\left(x\in\mathbb{S}\right)\mathcal{A}f\left(x\right)+I\left(x\in\partial\mathbb{S}\right)\mathcal{K}f\left(x\right).
	\end{equation}
\end{theorem}
The condition \eqref{eq:Wentzell} is known as the Wentzell boundary condition. The proof of this theorem is given in the appendix.

\section{CTMC Approximation}
\label{sec:simulation_2d_2_sticky_dimensions}

We will show how to construct a multidimensional CTMC to approximate a diffusion with stickiness in all dimensions and construction for the case where only some of the dimensions are sticky follows easily. To simplify the discussion, we set $\mathbb{S}=\{x\in\mathbb{R}^d:x^i>0,\ \textrm{for all}\ i=1,\ldots,\hat{d}\}$ and $\partial\mathbb{S}=\mathbb{R}^d_{+}\setminus\mathbb{S}$.  In the following, we develop two ways to construct the CTMC. 

\subsection{The Finite Difference Approach}
\label{subsec:simulation_2d_2_sticky_dimensions_algorithm}

The derivation of the transition rates for the CTMC is based on a discretization of the infinitesimal generator $\mathcal{G}$ given by \eqref{eq:generator}.
First, consider $x\in\mathbb{S}$. Let $e_i=(0,\ldots,0,1,0,\ldots,0)^\top\in\mathbb{R}^d$, i.e. the direction of the $i$-th coordinate, which is a possible direction for the Markov chain to move along. For simplicity, we use the same step size $h$ for all the directions but different step sizes can be used in our method with minor adjustments. 

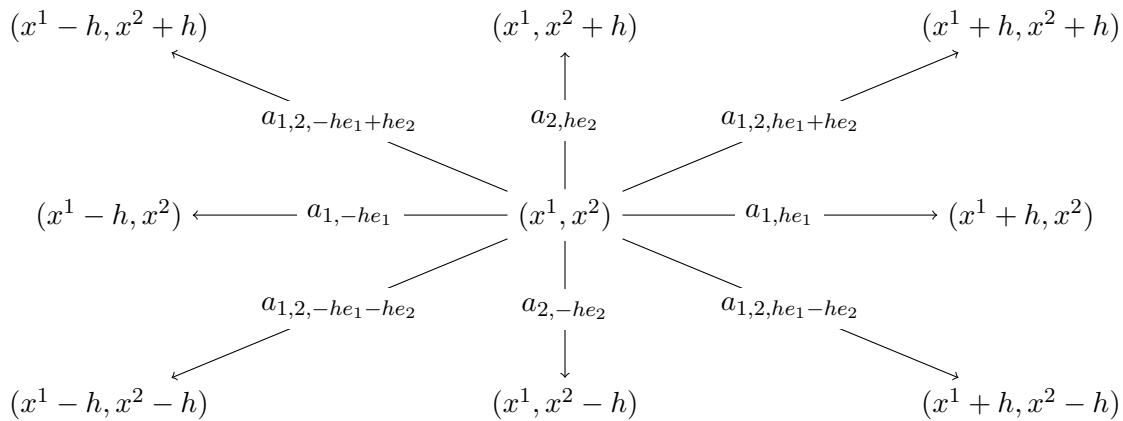
\begin{figure}[H]
	\centering
	\begin{tikzpicture}[node distance=2.5cm]
		\node [rectangle] (amp) {$(x^1-h,x^2+h)$};
		\node [rectangle, right of=amp, node distance=6cm] (anp) {$(x^1,x^2+h)$};
		\node [rectangle, right of=anp, node distance=6cm] (app) {$(x^1+h,x^2+h)$};
		\node [rectangle, below of=amp] (amn) {$(x^1-h,x^2)$};
		\node [rectangle, below of=anp] (ann) {$(x^1,x^2)$};
		\node [rectangle, below of=app] (apn) {$(x^1+h,x^2)$};
		\node [rectangle, below of=amn] (amm) {$(x^1-h,x^2-h)$};
		\node [rectangle, below of=ann] (anm) {$(x^1,x^2-h)$};
		\node [rectangle, below of=apn] (apm) {$(x^1+h,x^2-h)$};
		
		\draw[->] (ann) -- node [rectangle, fill=white] {$a_{1,2,-he_1+he_2}$} (amp);
		\draw[->] (ann) -- node [rectangle, fill=white] {$a_{2,he_2}$} (anp);
		\draw[->] (ann) -- node [rectangle, fill=white] {$a_{1,2,he_1+he_2}$} (app);
		\draw[->] (ann) -- node [rectangle, fill=white] {$a_{1,-he_1}$} (amn);
		\draw[->] (ann) -- node [rectangle, fill=white] {$a_{1,he_1}$} (apn);
		\draw[->] (ann) -- node [rectangle, fill=white] {$a_{1,2,-he_1-he_2}$} (amm);
		\draw[->] (ann) -- node [rectangle, fill=white] {$a_{2,-he_2}$} (anm);
		\draw[->] (ann) -- node [rectangle, fill=white] {$a_{1,2,he_1-he_2}$} (apm);
	\end{tikzpicture}
	\caption{Transition directions in two dimensions.} \label{fig:transitions_2D}
\end{figure}

We approximate partial derivatives like $\partial f/\partial x^i$ and $\partial^2 f/\partial (x^i)^2$ using central difference:
\begin{align*}
	\frac{\partial f}{\partial x^i}&=\frac{f(x+he_i)-f(x-he_i)}{2h},\\
	\frac{\partial^2 f}{\partial (x^i)^2}&=\frac{f\left(x+he_i\right)-2f\left(x\right)+f\left(x-he_i\right)}{h^2}.
\end{align*}
For the cross derivatives ($i\neq j$), we approximate them as
\begin{align*}
	\frac{\partial^2 f}{\partial x^i\partial x^j}&\approx\frac{2f\left(x\right)+f\left(x+he_i+he_j\right)+f\left(x-he_i-he_j\right)}{2h^2} \\
	&\qquad-\frac{f\left(x+he_i\right)+f\left(x-he_i\right)+f\left(x+he_j\right)+f\left(x-he_j\right)}{2h^2}\qquad\textrm{if}\ A^{i,j}\left(x\right)\geq 0, \\
	\frac{\partial^2 f}{\partial x^i\partial x^j}&\approx-\frac{2f\left(x\right)+f\left(x+he_i-he_j\right)+f\left(x-he_i+he_j\right)}{2h^2} \\
	&\qquad+\frac{f\left(x+he_i\right)+f\left(x-he_i\right)+f\left(x+he_j\right)+f\left(x-he_j\right)}{2h^2}\qquad\textrm{if}\ A^{i,j}\left(x\right)<0.
\end{align*}
The approximation depends on the sign of $A^{i,j}(x)$.
If it is positive, transitions along $e_i+e_j$ and $-e_i-e_j$ are more likely due to the positive correlation between dimensions $i$ and $j$. For the same reason, $-e_i+e_j$ and $e_i-e_j$ are used if $A^{i,j}<0$. See Figure \ref{fig:transitions_2D} for an illustration in the two-dimensional case. 

Putting all finite difference approximations together and reorganizing the terms, we obtain 
\begin{align*}
	\mathcal{G}f\left(x\right)&\approx\sum_{i=1}^d\left(-\frac{\mu^i\left(x\right)}{2h}+\frac{A^{i,i}\left(x\right)}{2h^2}-\sum_{j=1,j\neq i}^d\frac{\left\vert A^{i,j}\left(x\right)\right\vert}{2h^2}\right)f\left(x-he_i\right) \nonumber\\
	&\qquad+\sum_{i=1}^d\left(-\frac{A^{i,i}\left(x\right)}{h^2}+\sum_{j=1,j\neq i}^d\frac{\left\vert A^{i,j}\left(x\right)\right\vert}{h^2}\right)f\left(x\right)\nonumber \\
	&\qquad+\sum_{i=1}^d\left(\frac{\mu^i\left(x\right)}{2h}+\frac{A^{i,i}\left(x\right)}{2h^2}-\sum_{j=1,j\neq i}^d\frac{\left\vert A^{i,j}\left(x\right)\right\vert}{2h^2}\right)f\left(x+he_i\right)\nonumber \\
	&\qquad+\sum_{i,j=1,j\neq i}^d\frac{\left(A^{i,j}\left(x\right)\right)^+}{4h^2}\left(f\left(x+he_i+he_j\right)+f\left(x-he_i-he_j\right)\right) \nonumber\\
	&\qquad+\sum_{i,j=1,j\neq i}^d\frac{\left(A^{i,j}\left(x\right)\right)^-}{4h^2}\left(f\left(x+he_i-he_j\right)+f\left(x-he_i+he_j\right)\right).
\end{align*}
From this expression, we identify the transition rates of the CTMC as below:
\begin{align}
	a_{i,-he_i}(x)&=-\frac{\mu^i\left(x\right)}{2h}+\frac{A^{i,i}\left(x\right)}{2h^2}-\sum_{j=1,j\neq i}^d\frac{\left\vert A^{i,j}\left(x\right)\right\vert}{2h^2},\\
	a_{i,he_i}(x)&=\frac{\mu^i\left(x\right)}{2h}+\frac{A^{i,i}\left(x\right)}{2h^2}-\sum_{j=1,j\neq i}^d\frac{\left\vert A^{i,j}\left(x\right)\right\vert}{2h^2},\\
	a_{i,j,he_i+he_j}(x)&=a_{i,j,-he_i-he_j}(x)=\frac{\left(A^{i,j}\left(x\right)\right)^+}{2h^2},\\
	a_{i,j,he_i-he_j}(x)&=a_{i,j,-he_i+he_j}(x)=\frac{\left(A^{i,j}\left(x\right)\right)^-}{2h^2},
\end{align}
where $a_{i,\pm he_i}(x)$ is the rate of the CTMC transitioning to $x\pm he_i$ and $a_{i,j,\pm he_i\pm he_j}(x)$ is the rate transitioning to $x\pm he_i\pm he_j$.

The rates $a_{i,j,\pm he_i\pm he_j}(x)$ are always non-negative but $a_{i,\pm he_i}(x)$ may not be, leading to invalid Markov chains. In the following, we consider the case where the covariance matrix $A(x)$ is strictly diagonally dominant, i.e., 
for all $i=1,\ldots,d$ and $x\in\bar{\mathbb{S}}$,
\begin{equation*}
	A^{i,i}\left(x\right)-\sum_{j=1,j\neq i}^d\left\vert A^{i,j}\left(x\right)\right\vert> 0.
\end{equation*}
Under this condition, we have the following result. 
\begin{proposition}\label{prop:proper_transitions_diagonally_dominant}
Suppose $A$ is strictly diagonally dominant. For any given $x\in\bar{\mathbb{S}}$, the rates $a_{i,\pm he_i}(x)$ are non-negative if $h$ is small enough, and hence the Markov chain constructed by the finite difference approach is valid. 
\end{proposition}

Now consider the case where $x\in\partial\mathbb{S}$. Let $\mathcal{I}_0(x)=\{i\in\{1,\cdots,\hat{d}\}:x^i=0\}$, which is the set of sticky dimensions that are at zero and $\mathcal{I}_0(x)^c=\{1,\cdots,d\}\setminus\mathcal{I}_0(x)$. We discretize the generator as follows:
\begin{align*}
	\mathcal{G}f\left(x\right)&=\sum_{i=1}^d\hat{\beta}^i\left(x\right)\frac{\partial}{\partial x^i}f\left(x\right) +\frac{1}{2}\sum_{i,j=1}^{d}\hat{G}^{i,j}\left(x\right)\frac{\partial^2}{\partial x^i\partial x^j}f\left(x\right) \\
	&\approx\sum_{i\in\mathcal{I}_0(x)}\frac{\hat{\beta}^i\left(x\right)}{h}\left(f\left(x+he_i\right)-f\left(x\right)\right)\\
	&+\sum_{i\in\mathcal{I}_0(x)^c}\left(-\frac{\hat{\beta}^i\left(x\right)}{2h}+\frac{\hat{G}^{i,i}\left(x\right)}{2h^2}-\sum_{j\in\mathcal{I}_0(x)^c,j\neq i}\frac{\left\vert \hat{G}^{i,j}\left(x\right)\right\vert}{2h^2}\right)f\left(x-he_i\right) \nonumber\\
	&+\sum_{i\in\mathcal{I}_0(x)^c}\left(-\frac{\hat{G}^{i,i}\left(x\right)}{h^2}+\sum_{j\in\mathcal{I}_0(x)^c,j\neq i}\frac{\left\vert\hat{G}^{i,j}\left(x\right)\right\vert}{h^2}\right)f\left(x\right)\nonumber \\
	&+\sum_{i\in\mathcal{I}_0(x)^c}\left(\frac{\hat{\beta}^i\left(x\right)}{2h}+\frac{\hat{G}^{i,i}\left(x\right)}{2h^2}-\sum_{j\in\mathcal{I}_0(x)^c,j\neq i}\frac{\left\vert \hat{G}^{i,j}\left(x\right)\right\vert}{2h^2}\right)f\left(x+he_i\right)\nonumber \\
	&+\sum_{i,j\in\mathcal{I}_0(x)^c,j\neq i}\frac{\left(\hat{G}^{i,j}\left(x\right)\right)^+}{4h^2}\left(f\left(x+he_i+he_j\right)+f\left(x-he_i-he_j\right)\right) \nonumber\\
	&+\sum_{i,j\in\mathcal{I}_0(x)^c,j\neq i}\frac{\left(\hat{G}^{i,j}\left(x\right)\right)^-}{4h^2}\left(f\left(x+he_i-he_j\right)+f\left(x-he_i+he_j\right)\right).
\end{align*}
The resulting transition rates are given as follows:
	\begin{align*}
	a_{i,he_i}(x)&=\frac{\hat{\beta}^i\left(x\right)}{h},\qquad a_{i,-he_i}(x)=0, \qquad\qquad\qquad\qquad i\in\mathcal{I}_0(x), \\
	a_{i,\pm he_i}(x)&=\pm\frac{\hat{\beta}^i\left(x\right)}{2h}+\frac{\hat{G}^{i,i}\left(x\right)}{2h^2}-\sum_{j\in\mathcal{I}_0(x)^c,j\neq i}\frac{\left\vert\hat{G}^{i,j}\left(x\right)\right\vert}{2h^2}, \qquad i\in\mathcal{I}_0(x)^c \\
	a_{i,j,he_i+he_j}(x)&=a_{i,j,-he_i-he_j}(x)=\frac{\left(\hat{G}^{i,j}\left(x\right)\right)^+}{2h^2},\quad\qquad\qquad i,j\in\mathcal{I}_0(x)^c,\ j\neq i, \\
	a_{i,j,he_i-he_j}(x)&=a_{i,j,-he_i+he_j}(x)=\frac{\left(\hat{G}^{i,j}\left(x\right)\right)^-}{2h^2},\qquad\qquad\ i,j\in\mathcal{I}_0(x)^c,\ j\neq i, \\
	a_{i,j,\pm he_i\pm he_j}(x)&=0, \qquad\qquad\qquad\qquad\qquad\qquad\qquad\qquad\quad i\in\mathcal{I}_0(x),\ j\neq i.
\end{align*}
The rates $a_{i,\pm he_i}(x)$ are nonnegative with $h$ small enough if $\hat{G}$ is strictly diagonally dominant.

\subsection{The Eigendecomposition Approach}
\label{subsec:eigendecomposition_approach}

In the finite difference approach, at most two coordinates of the Markov chain can move at each transition. When several coordinates have strong correlations, the covariance matrix for the diffusion part is not strictly diagonally dominant, and hence the finite difference approach cannot guarantee valid rates. In such case, the co-movement in more than two coordinates should be allowed for the Markov chain to capture the strong correlations. 

We propose an alternative approach in which some of the moving directions are identified according to the covariance matrix. They are given by the eigenvectors of this matrix, which are orthogonal to each other. Specifically, as $A$ is positive definite, it can be written as $A=U\Lambda U^\top$, where $U$ is the matrix of normalized eigenvectors and $\Lambda$ a diagonal matrix containing the eigenvalues $\lambda_1,\ldots,\lambda_d$ which are all positive. In general, $A$ depends on $x$ and so do $U$ and $\Lambda$, but for simplicity this dependence is suppressed in the notations.  We can further rewrite $A$ as
\begin{equation}\label{eq:A-decomp}
	A=\sum_{i=1}^d\lambda_iu_iu_i^\top,
\end{equation}
where $u_i$ is the normalized eigenvector associated with $\lambda_i$.
From the results of principal component analysis, the eigenvector associated with the largest eigenvalue gives the direction along which the process varies the most, and the one associated with the second largest eigenvalue gives the second most variable direction, etc. Thus, the eigenvectors are natural directions to move the Markov chain.  

We first introduce a general idea. We define for each point $x$ a set of transition directions $M(x)=\{v_i(x):i\leq m(x)\}$, where $m(x)$ is uniformly bounded in $x$.
The transition rates $a_{i,hv_i}^d(x)$ and $a_{i,hv_i}^n(x)$ are used for matching the behavior of the drift term and the diffusion term, respectively.
These rates need to satisfy the following set of equations for $x\in\mathbb{S}$:
\begin{align}
	\sum_{i\in M(x)}a_{i,hv_i}^d(x)hv_i(x)&=\mu(x), \\
	\sum_{i\in M(x)}a_{i,hv_i}^n(x)\left(hv_i(x)\right)\left(hv_i(x)\right)^\top&=A(x),\\
	\sum_{i\in M(x)}a_{i,hv_i}^n(x)hv_i(x)&=0.\label{eq:condition_noise_rates}
\end{align}
The last equation holds because the first moment of the diffusion term is zero. For $x\in\partial\mathbb{S}$, we replace $\mu$ and $A$ with $\hat{\beta}$ and $\hat{G}$.

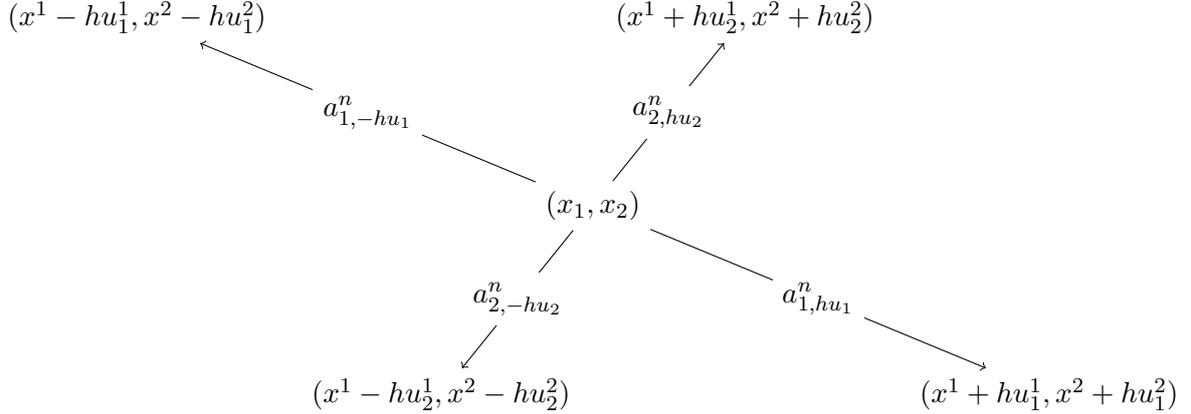
\begin{figure}[H]
	\centering
	\begin{tikzpicture}[node distance=2.5cm]
		\node [rectangle] (amp) {$(x^1-hu_1^1,x^2-hu_1^2)$};
		\node [rectangle, right of=amp, node distance=6cm] (anp) {};
		\node [rectangle, right of=anp, node distance=2cm] (app) {$(x^1+hu_2^1,x^2+hu_2^2)$};
		\node [rectangle, below of=anp] (ann) {$(x_1,x_2)$};
		\node [rectangle, below of=ann] (anm) {};
		\node [rectangle, right of=anm, node distance=6cm] (apn) {$(x^1+hu_1^1,x^2+hu_1^2)$};
		\node [rectangle, left of=anm, node distance=2cm] (anp) {$(x^1-hu_2^1,x^2-hu_2^2)$};
		
		\draw[->] (ann) -- node [rectangle, fill=white] {$a^n_{1,-hu_1}$} (amp);
		\draw[->] (ann) -- node [rectangle, fill=white] {$a^n_{2,hu_2}$} (app);
		\draw[->] (ann) -- node [rectangle, fill=white] {$a^n_{1,hu_1}$} (apn);
		\draw[->] (ann) -- node [rectangle, fill=white] {$a^n_{2,-hu_2}$} (anp);
	\end{tikzpicture}
	\caption{Transition directions for the diffusion part in two dimensions.} \label{fig:transitions_2D_noise_part}
\end{figure}

All possible solutions to the previous set of equations yield a Markov chain given that the rates are non-negative.
In order to limit the set of possible solutions and ensure the validity of the rates, we use the following set of transition directions:
\begin{equation}
	M\left(x\right)=\begin{cases} \{\mu\left(x\right)\}\cup\{u_i(x),-u_i(x):i=1,\ldots,d\} &\qquad\textrm{if}\ x\in\mathbb{S}, \\  \{\hat{\beta}\left(x\right)\}\cup\{u_i^{\hat{G}}(x),-u_i^{\hat{G}}(x):i=1,\ldots,d\} &\qquad\textrm{if}\ x\in\partial\mathbb{S}. \end{cases}
\end{equation}
Hence, the possible transition directions are given by the drift vector and the normalized eigenvectors of the covariance matrix pointing in both positive and negative directions. We can then obtain the following simple solutions:
\begin{align}\begin{split}\label{eq:transition_rates_general_drift}
		a_{h\mu}^d(x)&=\frac{1}{h},\qquad a_{i,\pm hu_i}^d(x)=0\qquad\qquad\qquad\quad\ \ \, \textrm{for}\ i=1,\ldots,d, \\
		a_{h\mu}^n(x)&=0,\qquad a_{i,-hu_i}^n(x)=a_{i,hu_i}^n(x)=\frac{\lambda_i}{2h^2}\qquad \textrm{for}\ i=1,\ldots,d.
\end{split}\end{align}
Note that by setting the rates $a_{i,-hu_i}^n(x)=a_{i,hu_i}^n(x)$, the condition in \eqref{eq:condition_noise_rates} is satisfied, because
\begin{equation}
	\sum_{i\in M(x)}a_{i,hv_i}^n(x)hv_i(x)=a_{h\mu}^n(x)+\sum_{i=1}^d\left(a_{i,-hu_i}^n(x)\left(-hu_i\right)+a_{i,hu_i}^n(x)\left(hu_i\right)\right)=0.
\end{equation}
An illustration of the possible movements of the chain capturing the diffusion term in the two-dimensional case is given in Figure \ref{fig:transitions_2D_noise_part}.

Recall that for $x\in\partial\mathbb{S}$, the $i$-th row and $i$-th column of the matrix $\hat{G}$ are equal to zero for all $i\in\mathcal{I}_0(x)$.
Hence, the eigenvalues of this diffusion matrix are $\lambda_i^{\hat{G}}>0$ ($i\in\mathcal{I}_0(x)^c$) and $\lambda_i^{\hat{G}}=0$ ($i\in\mathcal{I}_0(x)$).
We denote the resulting eigenvectors by $u_i^{\hat{G}}$ and let $u_i^{\hat{G}}=e_i$ for all $i\in\mathcal{I}_0(x)^c$.
The transition rates are then given as:
\begin{align}\begin{split}\label{eq:transition_rates_general_drift_boundary}
	a_{h\hat{\beta}}^d(x)&=\frac{1}{h},\qquad a_{i,\pm hu_i^{\hat{G}}}^d(x)=0\qquad\qquad\qquad\quad\ \ \, \textrm{for}\ i=1,\ldots,d \\
	a_{h\hat{\beta}}^n(x)&=0,\qquad a_{i,-hu_i^{\hat{G}}}^n(x)=a_{i,hu_i^{\hat{G}}}^n(x)=\frac{\lambda_i^{\hat{G}}}{2h^2}\qquad \textrm{for}\ i=1,\ldots,d.
\end{split}\end{align}

The rates defined in \eqref{eq:transition_rates_general_drift} and \eqref{eq:transition_rates_general_drift_boundary} are all non-negative independently of the structure of $A$ and $\hat{G}$.
The rate $a_{h\mu}^d(x)$ also does not depend on the sign of $\mu$ at the point $x$ and it captures the movement of the drift term accurately.
This also holds true for $\hat{\beta}$ whenever $x$ is at the boundary.

\begin{figure}[htbp!]
	\centering
	\begin{subfigure}{.496\textwidth}
		\centering
		\includegraphics[width=\linewidth]{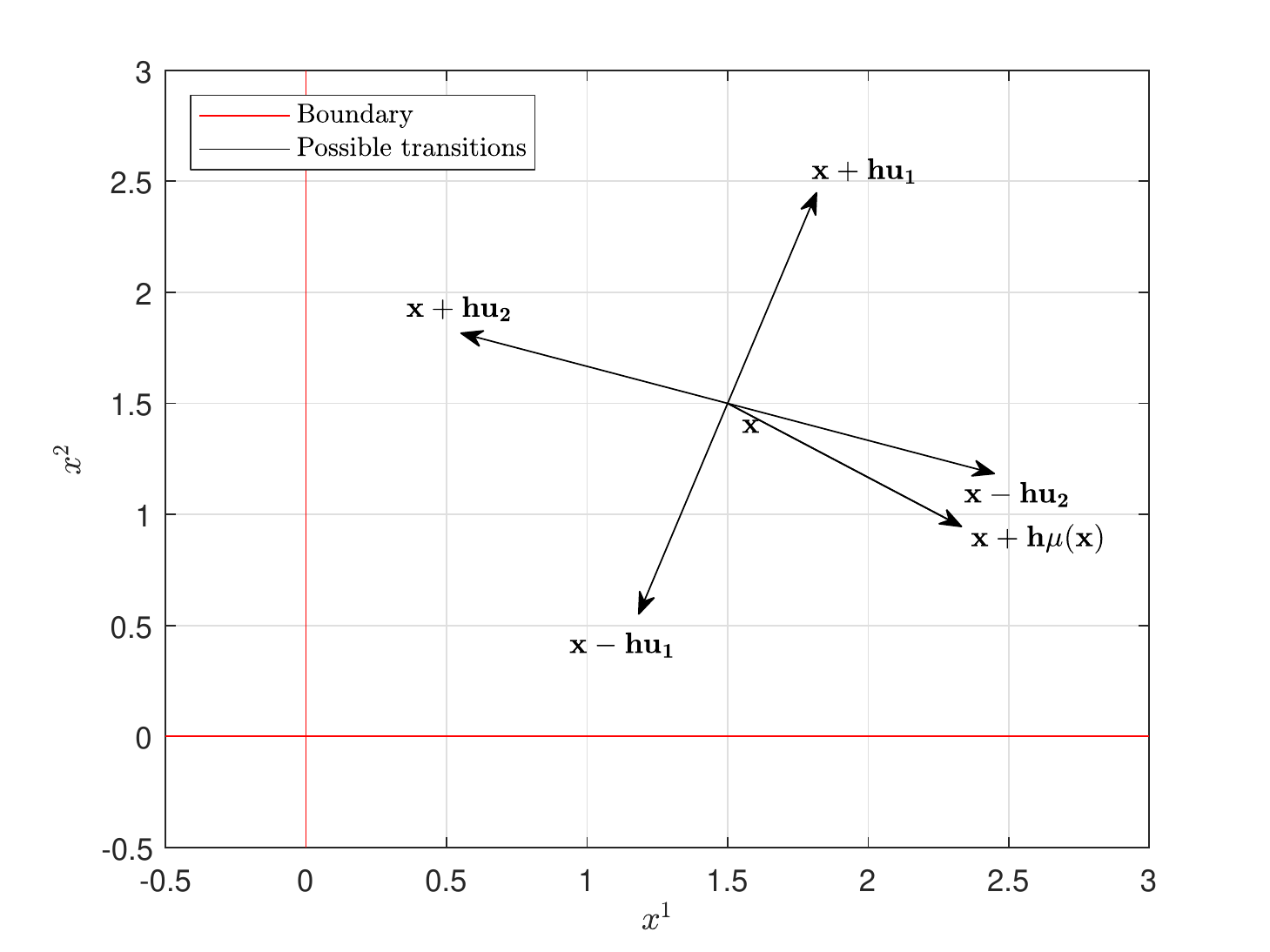}
	\end{subfigure}
	\begin{subfigure}{.496\textwidth}
		\centering
		\includegraphics[width=\linewidth]{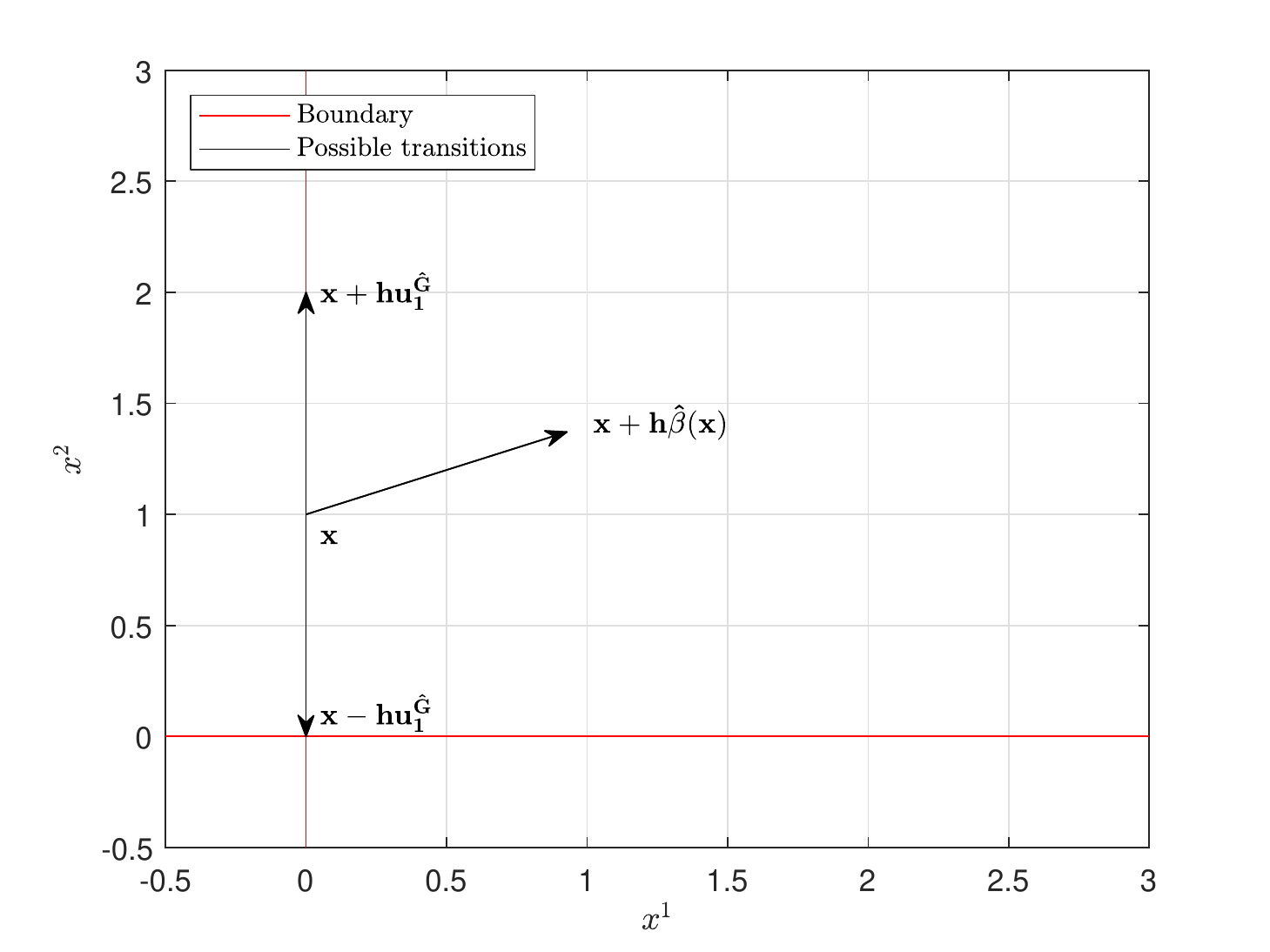}
	\end{subfigure}
	\begin{subfigure}{.496\textwidth}
		\centering
		\includegraphics[width=\linewidth]{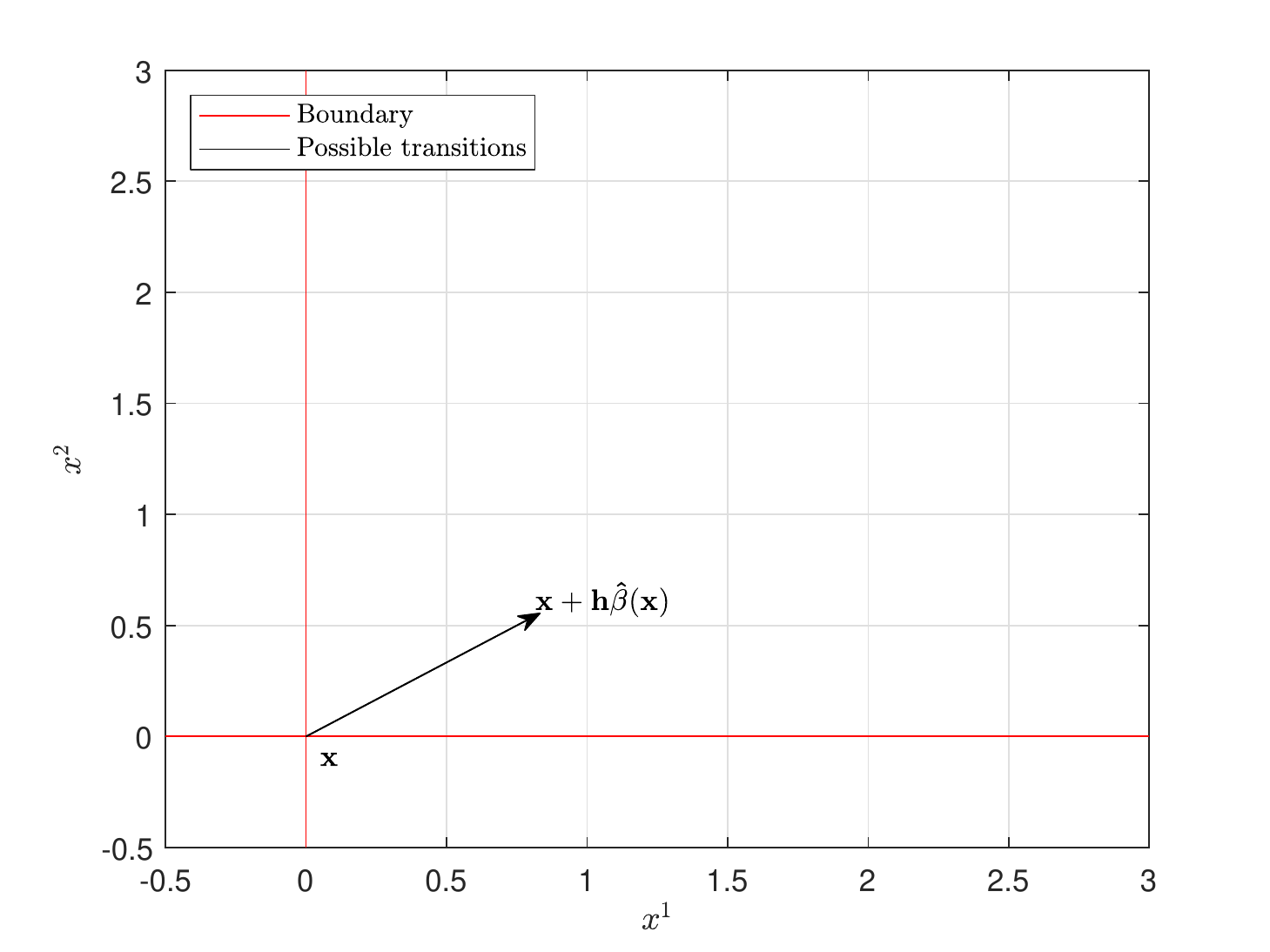}
	\end{subfigure}
	\begin{subfigure}{.496\textwidth}
		\centering
		\includegraphics[width=\linewidth]{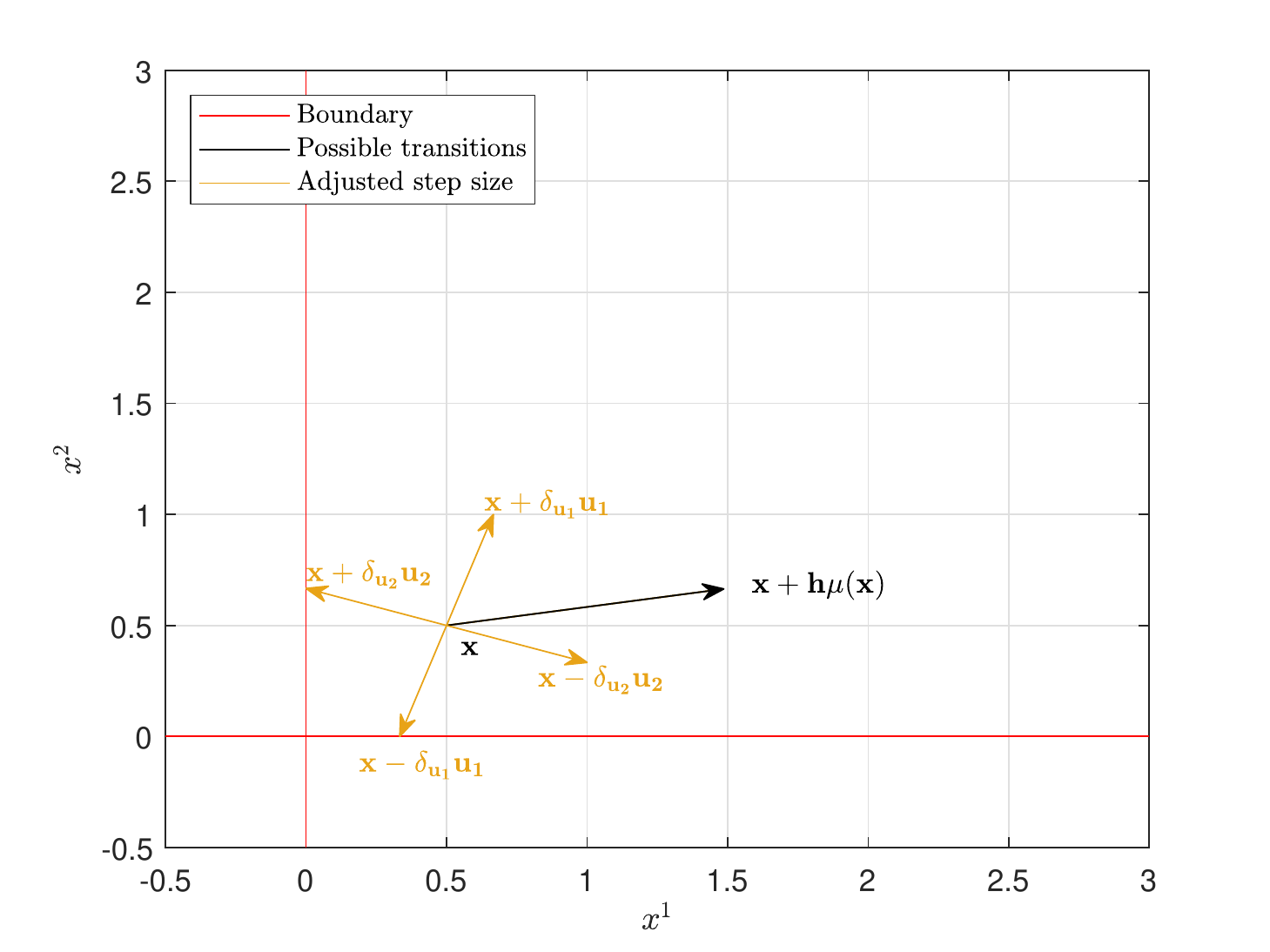}
	\end{subfigure}
	\caption{Possible transitions for an interior point (top left), a boundary point (top right), the origin (bottom left) and an interior point close to the boundary (bottom right).}
	\label{fig:2D_queuing_model_eigendecomposition_transition_rates}
\end{figure}

See Figure \ref{fig:2D_queuing_model_eigendecomposition_transition_rates} for the visualization of possible transitions of interior points and boundary points in two dimensions. It is worth noting that in the 2D problem, the eigendecomposition of $\hat{G}$ becomes trivial and no calculations are needed. For example, in the top right plot of Figure \ref{fig:2D_queuing_model_eigendecomposition_transition_rates}, we have $x^1=0$ and $x^2>0$. In this case, it is obvious that $\lambda_1^{\hat{G}}=0$ and $\lambda_2^{\hat{G}}=\hat{G}_{2,2}$. The normalized eigenvector associated with $\lambda_2^{\hat{G}}$ is simply $e_2$. In the bottom left plot of Figure \ref{fig:2D_queuing_model_eigendecomposition_transition_rates}, $x$ is the origin. In this case $\hat{G}(x)=0$, and hence the eigenvectors are all zero vectors. 

\subsection{Adjustment of Step Sizes}

It is important to note that if $x\in\mathbb{S}$ is close to the boundary, then a move along a direction with step size $h$ may lead the Markov chain out of $\bar{\mathbb{S}}$, and hence the step size must be trimmed. See the bottom right plot in Figure \ref{fig:2D_queuing_model_eigendecomposition_transition_rates} for an illustration. To adjust the step size along a direction, we must know the distance of the current point to the boundary $\partial\mathbb{S}$ along this direction.

In general,  the boundary is given by $\partial\mathbb{S}=\{x\in\mathbb{R}^d:\Phi\left(x\right)=0\}$. At an interior point $x$, for any chosen transition direction $u$, one computes the distance of $x$ to $\partial\mathbb{S}$ along $u$ (denoted by $\delta$) by solving $\Phi(x+\delta u)=0$, which is given by the smallest positive root (denoted by $\delta^u_{\text{min}}$).
The adjusted step size along $u$ is taken as the minimum of $h$ and $\delta^u_{\text{min}}$. 

When $\bar{\mathbb{S}}=\{x\in\mathbb{R}^d:x^1,\ldots,x^{\hat{d}}\geq 0\}$, $\Phi(x)$ is given by \eqref{eq:definition_Phi_orthant} and we can find the roots easily. In this case, we have
\begin{equation}\label{eq:definition_delta_x}
	\delta^u_{\text{min}}=\min\left\{-\frac{x^i}{u^i}: -\frac{x^i}{u^i}>0, i=1,\cdots,\hat{d}\right\}.
\end{equation}

After the step size of a direction is adjusted, we need to ensure the local consistency condition still holds. In our implementation, for simplicity we do the following:
\begin{itemize}
	\item For the finite difference approach, we use the same step size for all directions, which is given by the minimum of the adjusted step sizes along all possible directions.   
	
	\item For the eigendecomposition approach, if $u$ is an eigenvector, the adjusted step sizes for $u$ and $-u$ are set as the same. However, the adjusted step sizes could be different for different eigendirections and the drift direction.   
\end{itemize}

\subsection{Comparison of the Two Approaches}
We compare the finite difference (FD) approach and the eigendecomposition (ED) approach in various aspects below. 

\begin{itemize} 
\item Validity of the transition rates: the ED approach always guarantees valid transition rates, whereas the FD approach cannot for problems with strong correlations.  

\item The number of transitions in a unit time interval: in both approaches, the total rate of moving out of a state is $O(d/h^2)$. Thus, the number of transitions in a unit time interval is $O(d/h^2)$ for both approaches.

\item The time complexity of simulating one transition: the complexity is proportional to the number of transition directions. For the FD approach, it is equal to $2(d+{d\choose 2})=d^2+d$ and for the ED approach, it is given by $2d+1$. 

\item Additional calculations: the ED approach requires computing the eigendecomposition of the covariance matrix and this needs to be done multiple times if it is state-dependent. Furthermore, obtaining the distance to the boundary along a direction requires more calculations in the ED approach. For the FD approach, to adjust the step size we only need to consider the coordinate directions along which the distance to the boundary is directly given by a coordinate. In contrast, in general an equation needs to be solved in the ED approach as the transition directions may not be coordinate ones. 

\item Accuracy: for the same $h$, the ED approach is more accurate than the FD approach because it uses directions along which the process varies the most. This is confirmed by the numerical examples in Section \ref{sec:numerical_examples}.
\end{itemize} 

Considering the computational efficiency alone, which approach is faster when using the same $h$ depends on the complexity of simulating one transition and additional calculations required. For simulating one transition, the cost of the FD approach is close to the ED approach for $d=2$ but could be much higher for $d>2$.
However, for $d=2$ the ED approach might be slower due to the additional calculations required.
But as $d$ becomes greater, the computational edge of FD becomes smaller and the ED approach may eventually be faster after $d$ gets big enough.
In Section \ref{sec:numerical_examples}, we show the performance of these two approaches in two-dimensional problems.

\subsection{Convergence Rate}
\label{subsec:simulation_2d_2_sticky_dimensions_convergence}

Let $X_t$ be the sticky diffusion given by \eqref{eq:multidimensional_sticky_SDE_adjusted} and $Y_t$ is a CTMC living on $\bar{\mathbb{S}}^h$ constructed from either the finite difference approach or the eigendecomposition approach with $Y_0=X_0$. 

To analyze the convergence rate, we utilize the result in \cite{zeng1994} which studies parabolic PDEs with Wentzell boundary condition and the semigroup theory (\cite{ethier2005}. To apply the results in \cite{zeng1994}, we use their setting to assume the state space $\bar{\mathbb{S}}_{\textrm{loc}}=[0,r]^d$ ($r>0$), which is bounded and all the dimensions are sticky at the zero boundary. We define the extended state space by including the time dimension as $\mathcal{S}=(0,T]\times\bar{\mathbb{S}}_{\textrm{loc}}$ for a terminal time $T$.
We consider the function space $C^{2,\alpha}(\mathcal{S})$ for $0<\alpha\leq 1$, which consists of functions $g$ such that $g,\,\partial_xg,\,\partial_{xx}g,\,\partial_tg$ are H\"older continuous with exponent $\alpha$. When $\alpha=1$, H\"older continuity becomes Lipschitz continuity. If a function $g(x)$ independent of $t$ defined on $\bar{\mathbb{S}}_{\textrm{loc}}$ is in $C^{2,\alpha}(\mathcal{S})$, we will write it as $g\in C^{2,\alpha}(\bar{\mathbb{S}}_{\textrm{loc}})$.

\begin{theorem}\label{th:convergence_rate_ctmc}
	Suppose Assumption \ref{assumptions:existence} holds, $\hat{\beta}(x)$ and $\hat{G}(x)$ are bounded and Lipschitz continuous on the boundary of $\bar{\mathbb{S}}_{\textrm{loc}}$, and the payoff function $f\in C^{2,1}(\bar{\mathbb{S}}_{\textrm{loc}})$ satisfying the Wentzell boundary condition. 
	For $(t,x)\in\mathcal{S}$, consider the value function 	$v\left(t,x\right)=\mathcal{P}_tf\left(x\right)$.
	Then we have $v\in C^{2,1}(\mathcal{S})$,
	and, for any $x\in\bar{\mathbb{S}}^h$, there holds
	\begin{equation}
		\left\vert\mathbb{E}_x\left(f\left(Y_T\right)\right)-\mathbb{E}_x\left(f\left(X_T\right)\right)\right\vert\leq Ch,
	\end{equation}
	where the constant $C>0$ is independent of $h$ and $x$.
\end{theorem}

\begin{remark}
We can also use the setting in \cite{tsapovska2008}, where the author considers $\bar{\mathbb{S}}=\{x\in\mathbb{R}^d: x_1\geq 0, x_i\in\mathbb{R}, i\neq 1\}$ and derives the property of the parabolic PDE with Wentzell boundary condition. In addition to Assumption \ref{assumptions:existence} and that $\hat{\beta}(x)$ and $\hat{G}(x)$ are bounded and Lipschitz continuous on the boundary, if we further assume  
	\begin{equation*}
	\sum_{i,j=1}^dA^{i,j}\left(x\right)\xi_i\xi_j\geq C_1\xi^\top\xi,\qquad\sum_{i,j=2}^d\hat{G}^{i,j}\left(y\right)\eta_i\eta_j\geq C_2\eta^\top\eta
\end{equation*}
for some $C_1,C_2>0$ for all $x\in\bar{\mathbb{S}}$, $y\in\partial\mathbb{S}$, $\xi\in\mathbb{R}^d$ and $\eta\in\mathbb{R}^{d-1}$, then we obtain 
$v\in C^{2,\alpha}((0,T]\times\bar{\mathbb{S}})$ for any $\alpha\in (0,1)$ from Theorem 1 in \cite{tsapovska2008} (note that Lipschitz continuity is assumed in Assumption \ref{assumptions:existence}). Consequently, we can prove that for any $\alpha\in(0,1)$, there holds
\begin{equation}
	\left\vert\mathbb{E}_x\left(f\left(Y_T\right)\right)-\mathbb{E}_x\left(f\left(X_T\right)\right)\right\vert\leq Ch^\alpha
\end{equation}
for some constant $C>0$ independent of $h$ and $x$. The result shows the convergence order is arbitrarily close to one. 
\end{remark}

Our approach offers a method to simulate the sticky diffusion by simulating the CTMC that approximates it, which can be done exactly without time discretization. However, bias is created by spatial discretization and it is first order in the discretization level $h$. This result shows that our scheme is comparable to the Euler scheme, whose bias is first order in the step size of time discretization.


\section{Two Applications}
\label{sec:numerical_examples}

We consider two applications in two dimensions to demonstrate the performance of our method. In particular, we will validate the convergence rate numerically and compare the computational efficiency of the finite difference approach and the eigendecomposition approach as well as exact and discrete time simulation of the CTMC.

\subsection{Queuing Systems with Exceptional Service Policy}
\label{subsec:simulation_2d_2_sticky_dimensions_queuing_model}
Our first application considers a multi-server queuing system where customers receive exceptional service whenever a server is idle. Studies on this type of system can be found in for example \cite{welch1964}, \cite{lemoine1974}, \cite{lemoine1975} and \cite{harrison1981}. Diffusions with sticky boundaries arise as the heavy traffic limit of such system (see \cite{racz2015}). Suppose there are $d$ servers in the system. The heavy traffic limit is given by
\begin{equation}\label{eq:sde-queue}	
	dX_t=\Sigma I\left(X_t\in\mathbb{S}\right)dB_{1,t}+\hat{\beta}\left(X_t\right)I\left(X_t\in\partial\mathbb{S}\right)dt.
\end{equation}
The drift at the boundary and the volatility in the interior are set as
\begin{alignat*}{2}
	\hat{\beta}(x)&=\sum_{i=1}^d\eta_iI\left(x^i=0\right) &\qquad\qquad\textrm{with}\ \eta_i=\left(\eta_i^1,\ldots,\eta_i^d\right)^\top\in\mathbb{R}_{>0}^d, \\ 
	\Sigma&=\sigma L &\qquad\qquad\textrm{for}\ x\in\bar{\mathbb{S}},\ \sigma>0,
\end{alignat*}
and $L$ is a $d$-by-$d$ lower triangular matrix, such that
\begin{equation}\label{eq:covariance_matrix_spacing_process}
	LL^\top=\begin{pmatrix}  2 & -1 & 0 & 0 & \cdots & 0 \\ -1 & 2 & -1 & 0 &\cdots & 0 \\ 0 & -1 & 2 & -1 & \cdots & 0 \\ \vdots & \vdots & \ddots & \ddots & \ddots & \vdots \\ 0 & \cdots & 0 & -1 & 2 & -1 \\ 0 & \cdots & 0 & 0 & -1 & 2 \end{pmatrix}.
\end{equation}
Hence, $X$ is a sticky Brownian motion in the positive orthant with stickiness in all dimensions.

\begin{remark}
	\cite{racz2015} considered a sticky $(d+1)$-dimensional diffusion and they mentioned that the spacing process of this sticky diffusion is the heavy traffic limit of a queuing system of $d$ servers with exceptional policy. Adapting their formulation to our setting leads to the form of the process given above. 
\end{remark}

To show a numerical example, we consider a two-server case ($d=2$) and set the parameters as follows: 
\begin{equation*}
	\eta=\begin{pmatrix} \eta_1^1 & \eta_2^1 \\ \eta_1^2 & \eta_2^2 \end{pmatrix} = \begin{pmatrix} 0.01 & 0.99 \\ 0.90 & 0.95 \end{pmatrix},\qquad\sigma=1.
\end{equation*}
These values are taken from Section 4 in \cite{doytchinov2001}.
We depart from the reference by making $\eta_1^1$ much smaller than the other entries to create more stickiness on the boundary $\{x\in\mathbb{R}^2: x^1=0\}$. To construct a CTMC approximation using the eigendecomposition approach, we calculate the eigenvalues and eigenvectors of the covariance matrix 
\begin{equation}
	A=\begin{pmatrix} 2 & -1 \\ -1 & 2 \end{pmatrix}, 
\end{equation}
which are given by
\begin{equation*}
	\lambda_1=1,\qquad \lambda_2=3,\qquad u_1=\begin{pmatrix} \frac{1}{\sqrt{2}} \\ \frac{1}{\sqrt{2}} \end{pmatrix},\qquad u_2=\begin{pmatrix} \frac{1}{\sqrt{2}} \\ -\frac{1}{\sqrt{2}} \end{pmatrix}.
\end{equation*}
Since $A$ is state-independent in this problem, the eigendecomposition only needs to be done once. 

We compute
\begin{equation*}
	v(T,x)=\mathbb{E}_x\left(f\left(X_T\right)\right)=\mathbb{E}_x\left(X_T^1+X_T^2\right),
\end{equation*}
for $x=(0,0)^\top$ and $T=1$. To provide an accurate benchmark for our simulation method, we use the alternate direction implicit (ADI) scheme (see e.g., Equation (20.9) in \cite{duffy2006}) to numerically solve the PDE with Wentzell boundary condition for $v(T,x)$. The resulting value is $v\left(T,x\right)=0.923377$, which is accurate for the first six decimal places.

\subsubsection{Sample Paths}
\label{subsec:simulation_2d_2_sticky_dimensions_queuing_model_results_fd}

\begin{figure}[htbp!]
	\centering
	\begin{subfigure}{.496\textwidth}
		\centering
		\includegraphics[width=\linewidth]{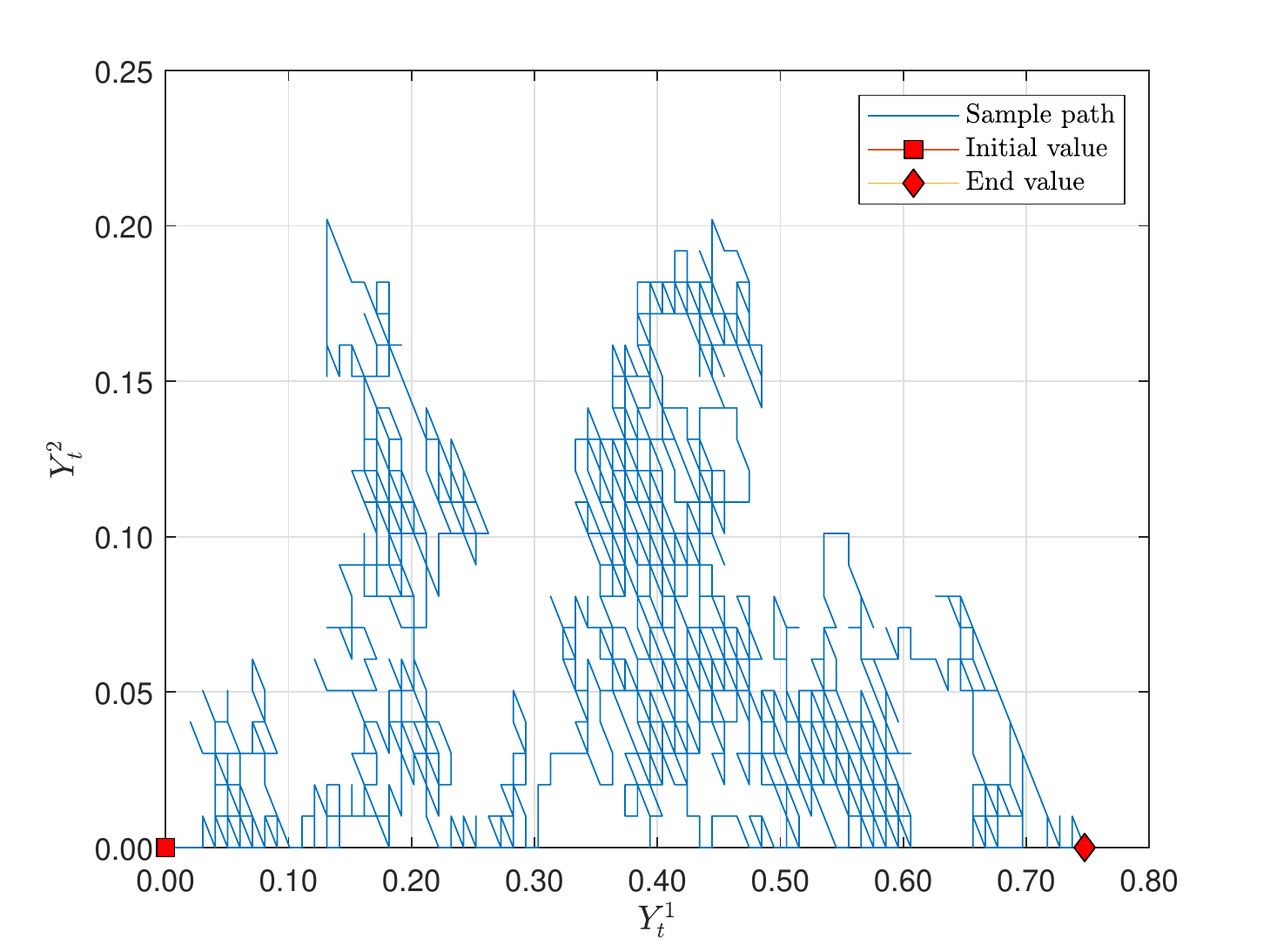}
	\end{subfigure}
	\begin{subfigure}{.496\textwidth}
		\centering
		\includegraphics[width=\linewidth]{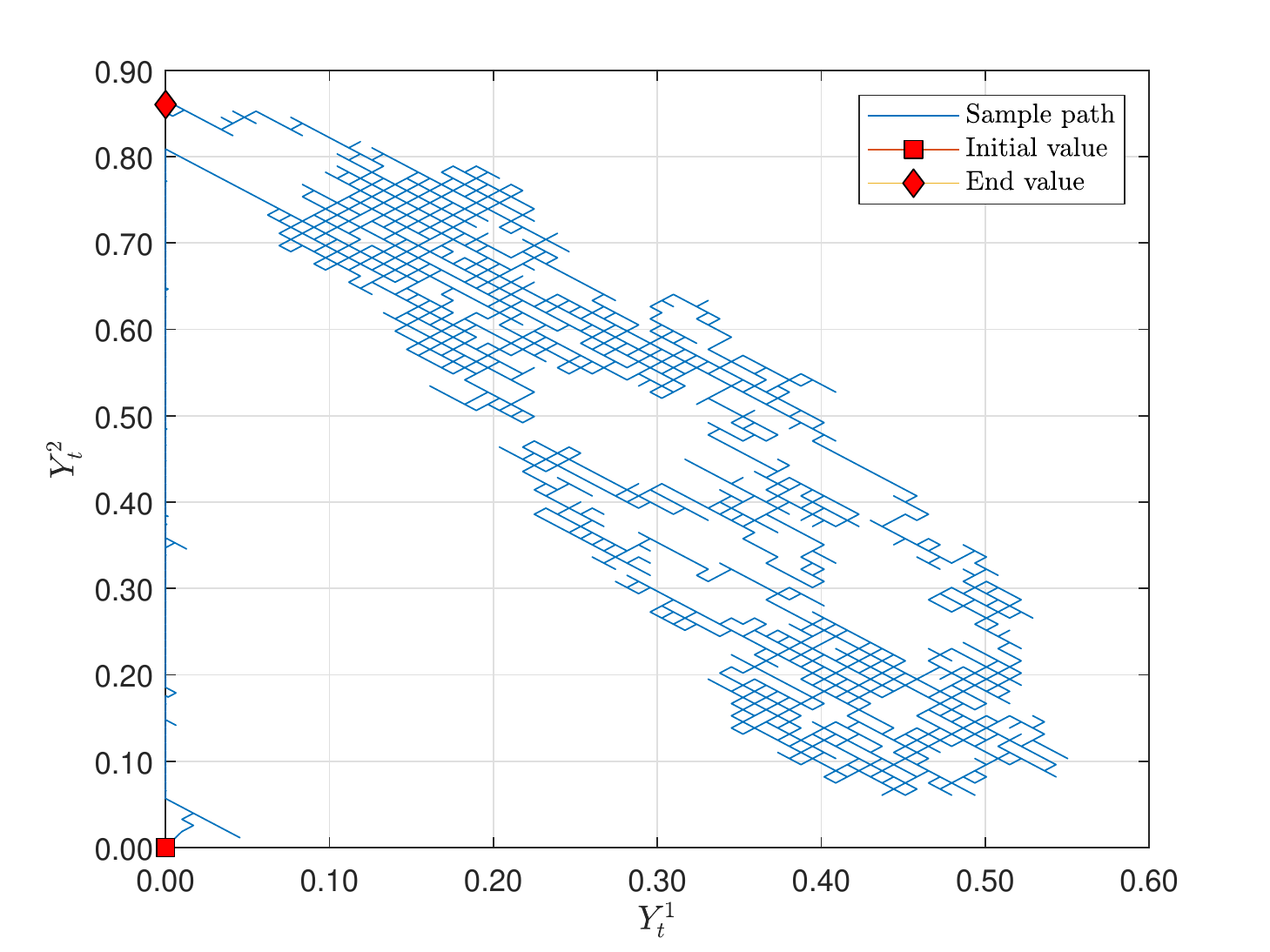}
	\end{subfigure}
	\caption{Sample paths over the time interval $[0,1]$ generated from the CTMC constructed by the finite difference approach (left) and the eigendecomposition approach (right) for the sticky Brownian motion \eqref{eq:sde-queue}. }
	\label{fig:2D_queuing_model_sample_paths}
\end{figure}

We first show paths generated by the finite difference approach and the eigendecomposition approach in Figure \ref{fig:2D_queuing_model_sample_paths}. In order to clearly show the differences, we choose a larger discretization of $h=1/100$. 

In the path generated by the finite difference approach, the Markov chain can move along $\pm e_1$, $\pm e_2$, $-e_1+e_2$, or $e_1-e_2$ (note that the correlation is negative in this model). In comparison, in the path generated by the eigendecomposition approach, the moving directions are $\pm u_1$ and $\pm u_2$ for interior points (movements along $\pm e_1$ and $\pm e_2$ are impossible because there is no drift in the interior) and they become $\pm e_1$ or $\pm e_2$ at the boundary. It is interesting to note that in the right plot, the process sticks to the boundary $\{x:x^1=0\}$ for a while once it hits there. This happens because $\eta^1_1$, the drift for the first coordinate to leave zero, is very small.

\begin{figure}[htbp!]
	\centering
	\begin{subfigure}{.496\textwidth}
		\centering
		\includegraphics[width=\linewidth]{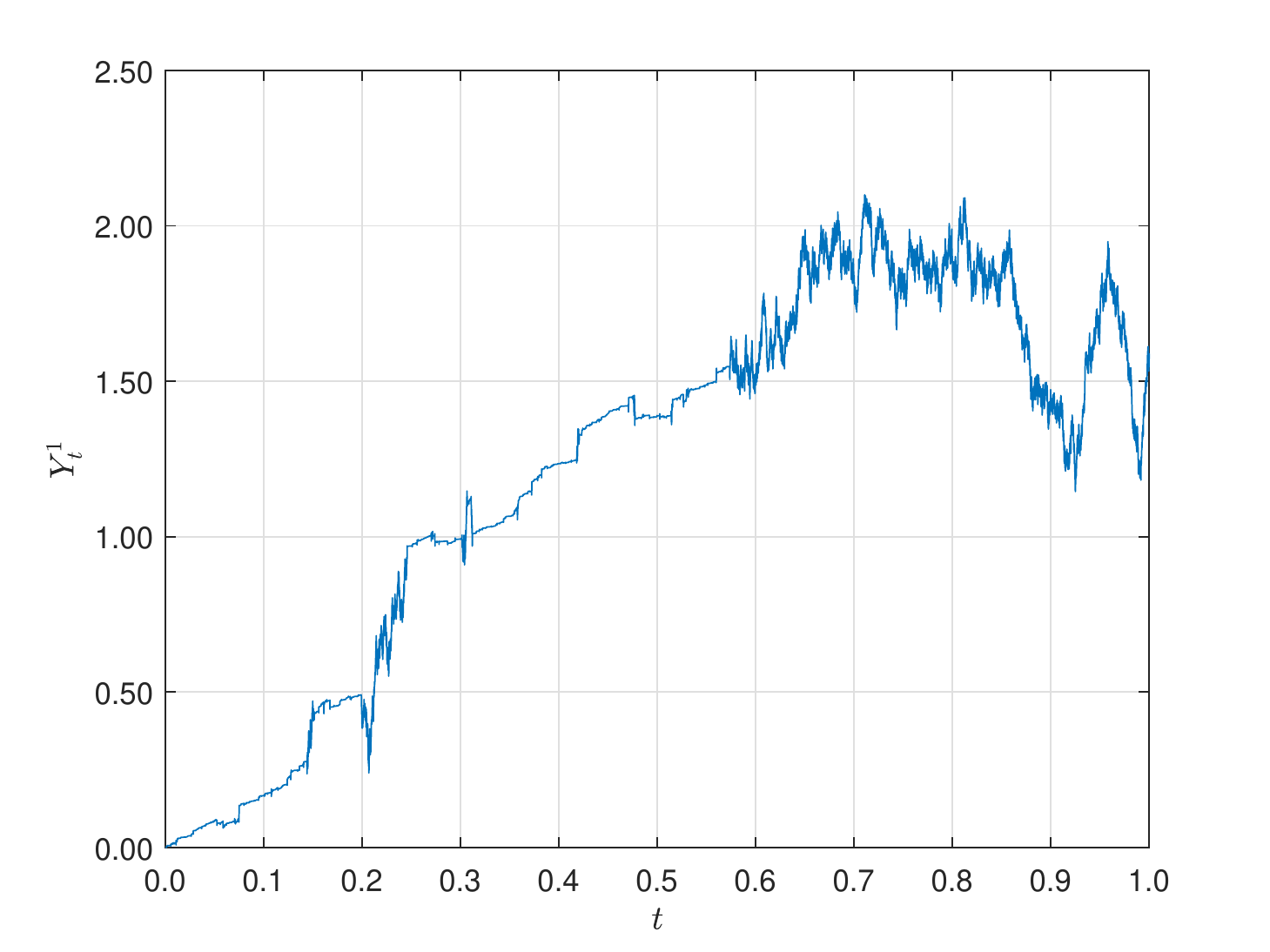}
	\end{subfigure}
	\begin{subfigure}{.496\textwidth}
		\centering
		\includegraphics[width=\linewidth]{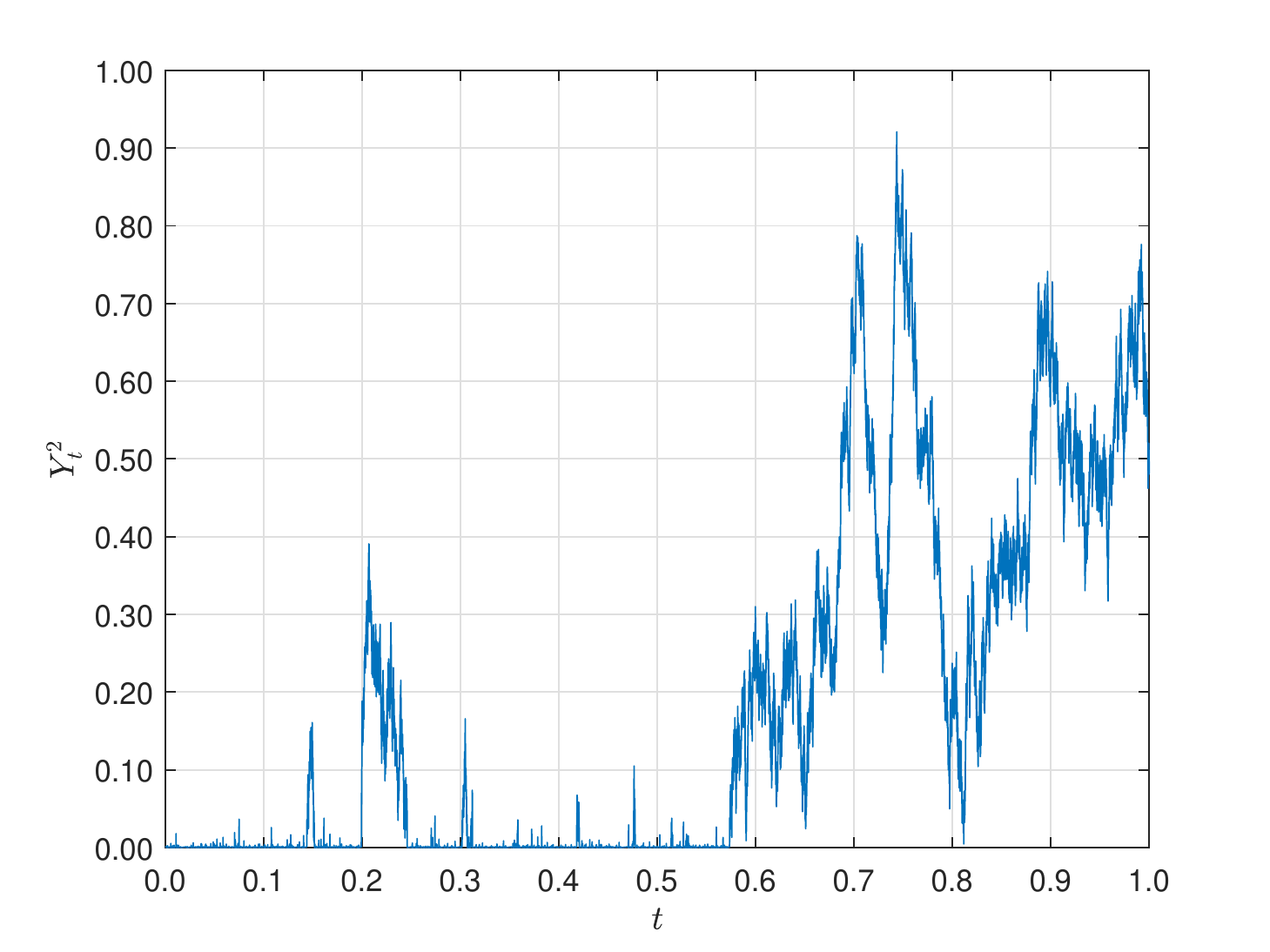}
	\end{subfigure}
	\begin{subfigure}{.496\textwidth}
		\centering
		\includegraphics[width=\linewidth]{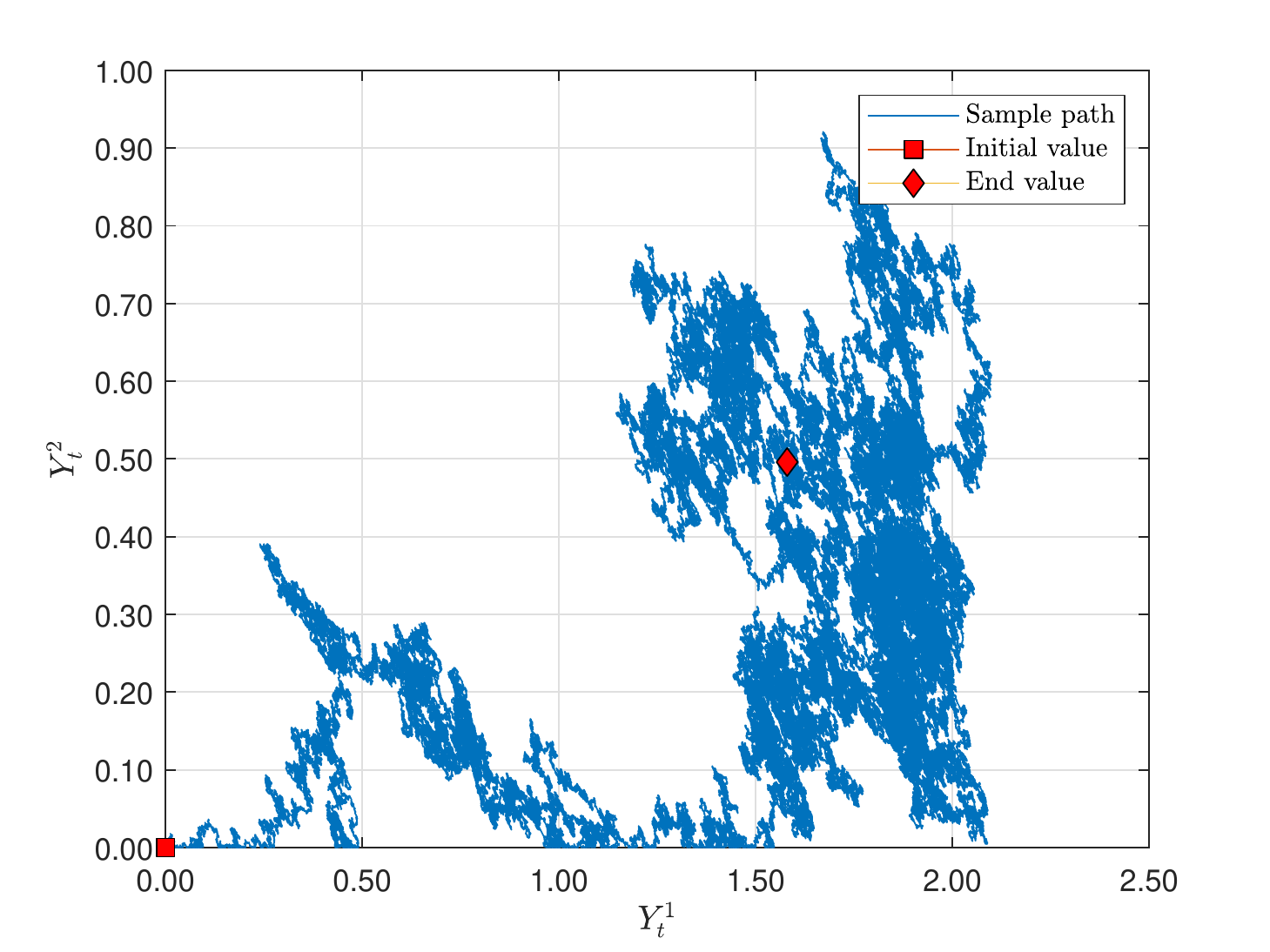}
	\end{subfigure}
	\caption{One sample path over the time interval $[0,1]$ generated from the CTMC with $h=1/1000$ constructed by the eigendecomposition approach for the sticky Brownian motion \eqref{eq:sde-queue}. The first row shows the two dimensions separately while the second row shows the movement of the process in the two dimensional state space. }
	\label{fig:2D_queuing_model_sample_paths_eigen}
\end{figure}

In Figure \ref{fig:2D_queuing_model_sample_paths_eigen}, we regenerate a path from the eigendecomposition approach with $h=1/1000$. This path provides a more accurate approximation to the true path of the two-dimensional sticky Brownian motion model given by \eqref{eq:sde-queue}.

\subsubsection{Convergence Results}
\label{subsubsec:convergence_results}

\begin{figure}[htbp!]
	\centering
	\begin{subfigure}{.496\textwidth}
		\centering
		\includegraphics[width=\linewidth]{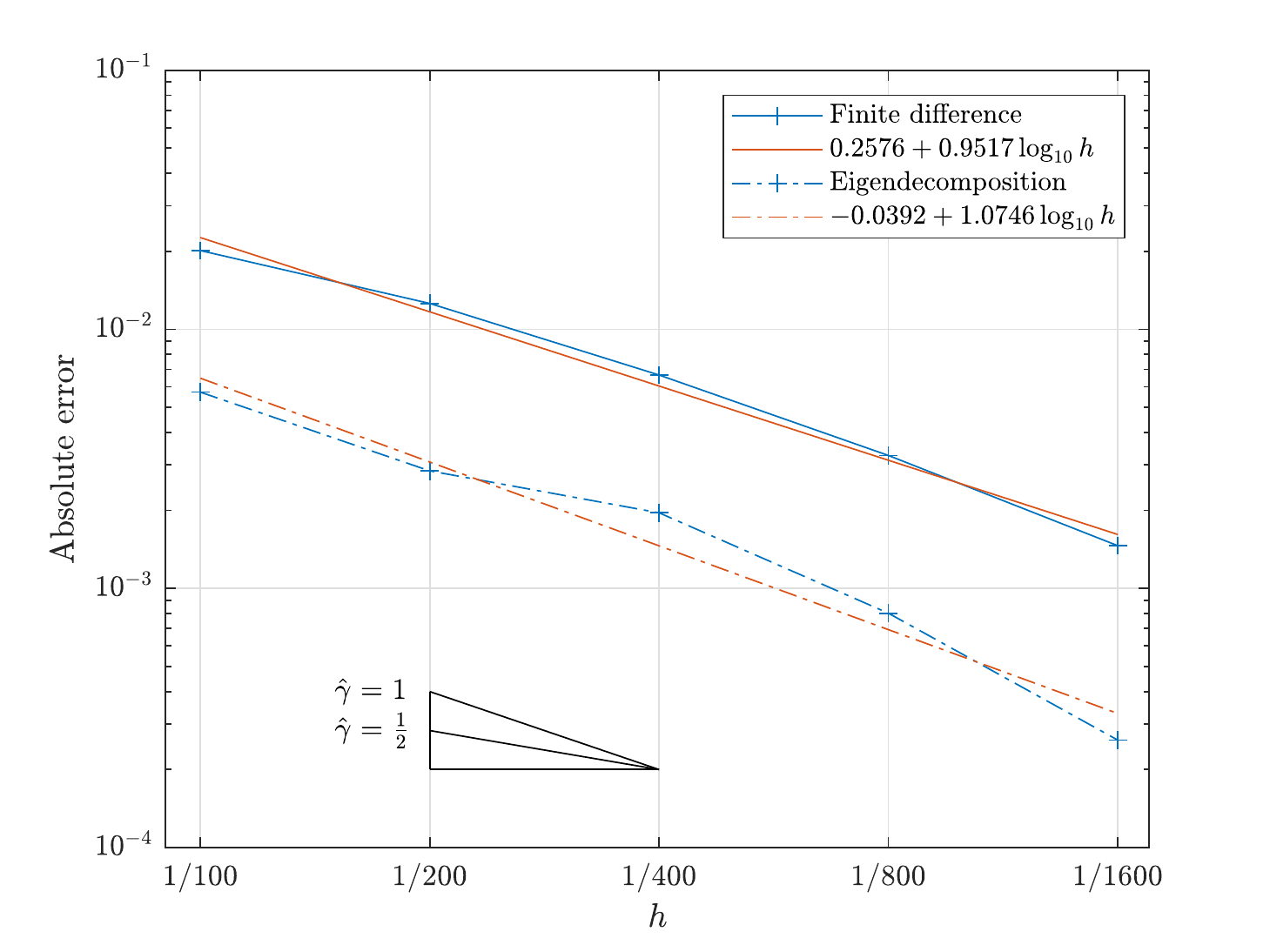}
	\end{subfigure}
	\begin{subfigure}{.496\textwidth}
		\centering
		\includegraphics[width=\linewidth]{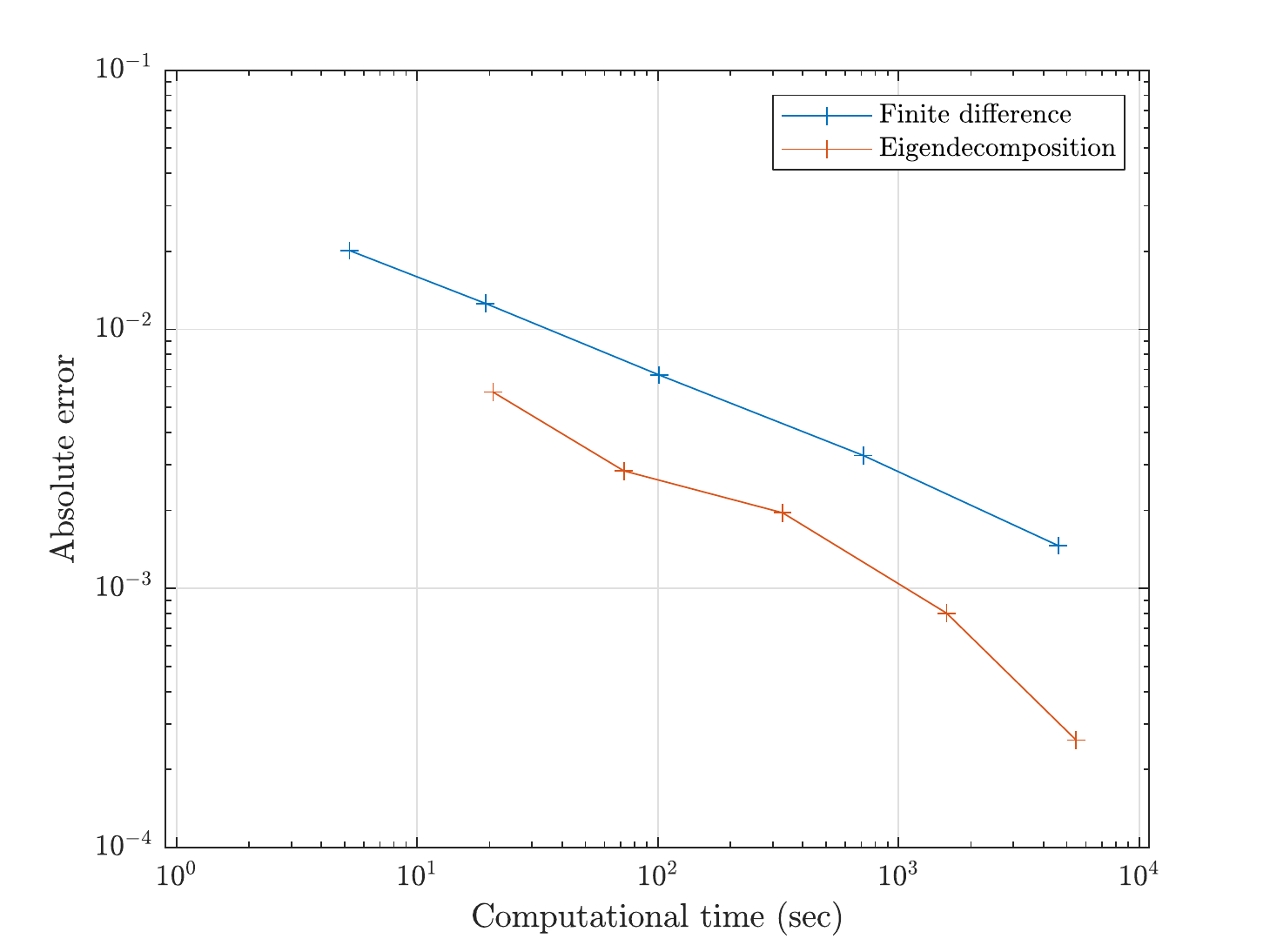}
	\end{subfigure}
	\caption{Convergence rate (left) and absolute error vs. computational time (right) for CTMC simulation of the sticky Brownian motion \eqref{eq:sde-queue}. Both plots are on log-log scale. }
	\label{fig:2D_queuing_model_convergence_results}
\end{figure}

To show convergence of our method, we set
\begin{equation}\label{eq:h}
	h=\frac{1}{100},\frac{1}{200},\frac{1}{400},\frac{1}{800},\frac{1}{1600},
\end{equation}
For each value of $h$, we generate $10^5$ paths from the CTMC to compute $v(T,x)$.
Figure \ref{fig:2D_queuing_model_convergence_results} shows the results for the finite difference and eigendecomposition approach in one plot to allow for easy comparison.
The numerically estimated convergence orders are 0.9517 for the finite difference approach and 1.0746 for the eigendecomposition approach, which are close to the theoretical convergence order of 1. The eigendecomposition approach has a smaller constant for the error, making it more accurate than the finite difference approach for the same level of $h$. The improvement in accuracy can be explained by the adjusted directions it uses for moving the CTMC.

The right plot in Figure \ref{fig:2D_queuing_model_convergence_results} shows the performance of these two approaches in terms of accuracy and computational time. For similar levels of accuracy the eigendecomposition approach is faster.
However, for a fixed level of $h$, the eigendecomposition approach is somewhat slower because additional computations related to adjusting the step size are performed.

\subsubsection{Exact vs. Discrete Time Simulation of the CTMC}
Although the CTMC can be simulated exactly, in general we need to simulate $O(d/h^2)$ number of transitions in a unit time interval, which can be too many. To save computations, we can simulate the CTMC in discrete time as follows. 

Consider a time grid $t_0=0, t_1=T/N, t_2=2T/N,\ldots,t_N=T$ for some terminal time $T$. Simulate 
$e\sim\mathrm{Exp}(1)$. Set $e_{t_0}=0$ and calculate $e_{t_{i+1}}=e_{t_i}+ a_0\tfrac{T}{N}$, where $a_0$ is the negative of the sum of transition rates out of the current state. The first transition time is given by  $\min\{t_i: e_{t_i}\geq e\}$, and then a state from $\mathbb{S}$ is drawn as the next state. Afterwards the same procedure starts over.

This discrete time scheme assumes that there is at most one transition in a time interval of length $T/N$, which inevitably introduces additional bias. The bias created by time discretization for each time step is $O(1/N^2)$ and the total bias is $O(1/N)$ over $[0,T]$. If we set $h=O(1/N)$, the overall bias from spatial and time discretization is still $O(h)$. Thus, simulating the CTMC in discrete time would not change the first order convergence of our method.

\begin{figure}[htbp!]
	\centering
	\begin{subfigure}{.496\textwidth}
		\centering
		\includegraphics[width=\linewidth]{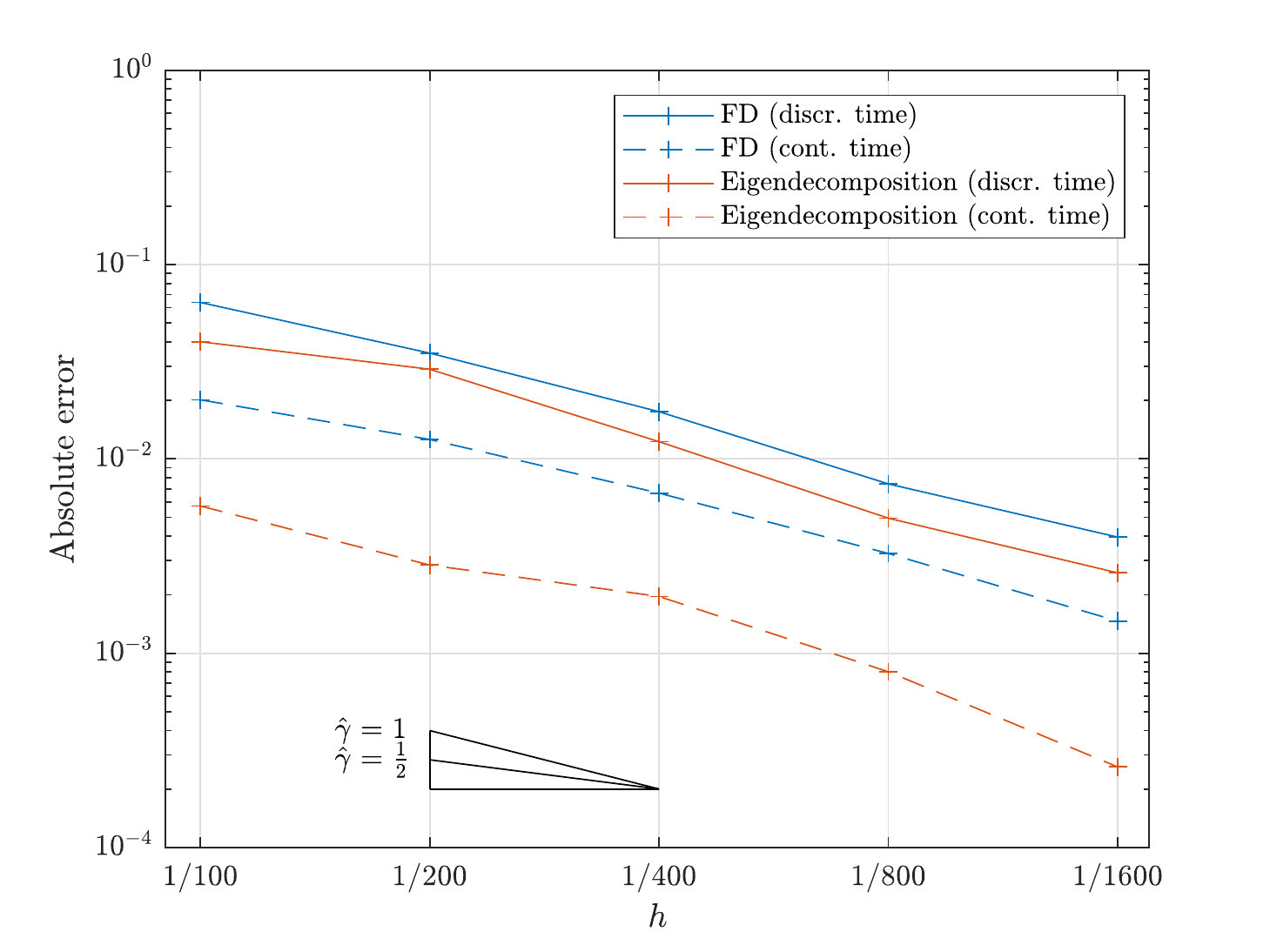}
	\end{subfigure}
	\begin{subfigure}{.496\textwidth}
		\centering
		\includegraphics[width=\linewidth]{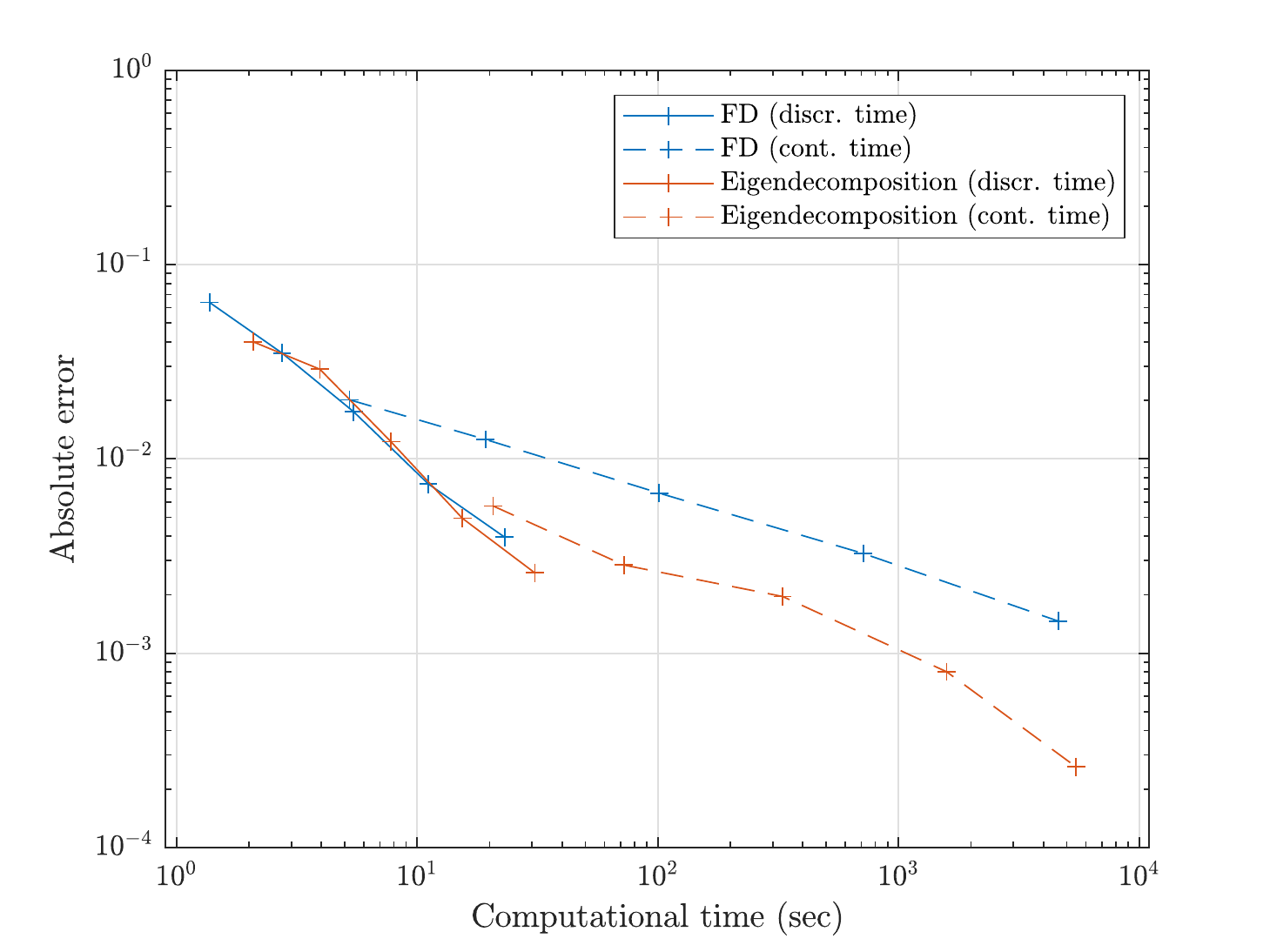}
	\end{subfigure}
	\caption{Convergence rate (left) and absolute error vs. computational time (right) for CTMC simulation of the sticky Brownian motion \eqref{eq:sde-queue}. Both plots are on log-log scale. Each marker in both plots corresponds to one level of $h$.}
	\label{fig:2D_queuing_model_convergence_results_discrete_time}
\end{figure}

Figure \ref{fig:2D_queuing_model_convergence_results_discrete_time} shows the results of exact and discrete time simulation of the CTMC using the two construction approaches. The discrete time simulation is less accurate for the same level of $h$ as expected. Nevertheless, from the right plot, we see that it can save computational time for obtaining similar levels of accuracy (although it cannot attain error levels of $10^{-3}$ using the values of $h$ considered), and this is because only $O(1/h)$ transitions are simulated over a unit time interval as opposed to $O(d/h^2)$ for the exact simulation. 

\subsection{A Multi-Factor Sticky Interest Rate Model}
\label{subsec:2D_interest_rate_model}
After the 2008 financial crisis, many central banks have been keeping a low-interest rate policy. This motivates the development of new interest rate models for a low interest environment. A natural way to model such phenomenon is creating stickiness for the short rate (instantaneous interest rate) at a low level such as zero, and the stickiness is determined by additional factors. This type of sticky short rate model is studied in detail in \cite{nie2020} and \cite{nie2017}, which explain why empirically it can be more realistic than the standard shadow rate model considered in \cite{kim2012}.

We consider the two-dimensional sticky short rate model in \cite{nie2017}. Let $X^1_t$ be the short rate at time $t$ and $X^2_t$ be a factor process. Set $\mathbb{S}=\mathbb{R}_{>0}\times\mathbb{R}$ and $\partial\mathbb{S}=\{0\}\times\mathbb{R}$. The model assumes 
\begin{align}\begin{split}\label{eq:2D_interest_rate_model}
		dX_t&=I\left(X_t\in\mathbb{S}\right)\left(K\left(\theta-X_t\right)dt+\Sigma dB_{1,t}\right) \\
		&+I\left(X_t\in\partial\mathbb{S}\right)\left(\begin{pmatrix} \nu\left(X^2_t\right) \\ \kappa_2\left(\theta_2-X^2_t\right)\end{pmatrix}dt+\begin{pmatrix} 0 & 0 \\ 0 & \sigma_2\end{pmatrix}dB_{2,t}\right),
\end{split}\end{align}
where $X_t=(X^1_t,X^2_t)\in\mathbb{R}^2$, $K,\Sigma\in\mathbb{R}^{2\times 2}$, $\theta,B_{1,t},B_{2,t}\in\mathbb{R}^2$ and $\kappa_2,\theta_2,\sigma_2\in\mathbb{R}$.
Here, $B_{1,t}$ and $B_{2,t}$ are two-dimensional standard Brownian motions which are independent of each other. In the interior of the state space, this model assumes $(X^1,X^2)$ follows a two-dimensional Ornstein-Uhlenbeck (OU) process. Once $X^1$ reaches zero, it is sticky there and the stickiness is inversely related to the function $\nu$, which is given by 
\begin{equation}
	\nu\left(X_t^2\right)=\frac{\nu}{1+\exp\left(-100X_t^2\right)}.
\end{equation}
The factor $X^2$ is not sticky and it is unbounded. But its dynamics differs for $X^1>0$ and $X^1=0$. \cite{nie2017} studies the empirical performance of this model in details. In practice, if one factor is not enough, one can extend the model by making the factor $X^2$ multidimensional. 

To develop a numerical example, we use the parameter values given in Section 7.6 of \cite{nie2017}, which are obtained by fitting the model to daily US yield curve data over a long period. They are given by
\begin{align}
	K&=\begin{pmatrix} 0.3076 & -0.1943 \\ -0.0401 & 0.0198\end{pmatrix},\qquad \theta=\begin{pmatrix} 0.0008 \\ -0.0363 \end{pmatrix},\qquad \Sigma=\begin{pmatrix} 0.0253 & 0 \\ 0 & 0.0189\end{pmatrix}, \\
	\kappa_2&=0.0665,\qquad \theta_2=0.0134,\qquad \sigma_2=0.1051,\qquad \nu=0.0079.
\end{align}
We compute the price of a zero-coupon bond with unit face value and maturity $T$, which is given by 
\begin{equation}
	v\left(T,x\right)=\mathbb{E}_x\left(\exp\left(-\int_0^{T}X^1_sds\right)\right).
\end{equation}
In our example, we set $T=1$ and $x=(0.01,0)^\top$. The benchmark is again computed by a finite difference scheme using ADI (see Section 7.2 and B.2 in \cite{nie2017}). The resulting value is 0.988626, which is accurate to the sixth decimal place.

\subsubsection{Sample Paths}
\begin{figure}[htbp!]
	\centering
	\begin{subfigure}{.496\textwidth}
		\centering
		\includegraphics[width=\linewidth]{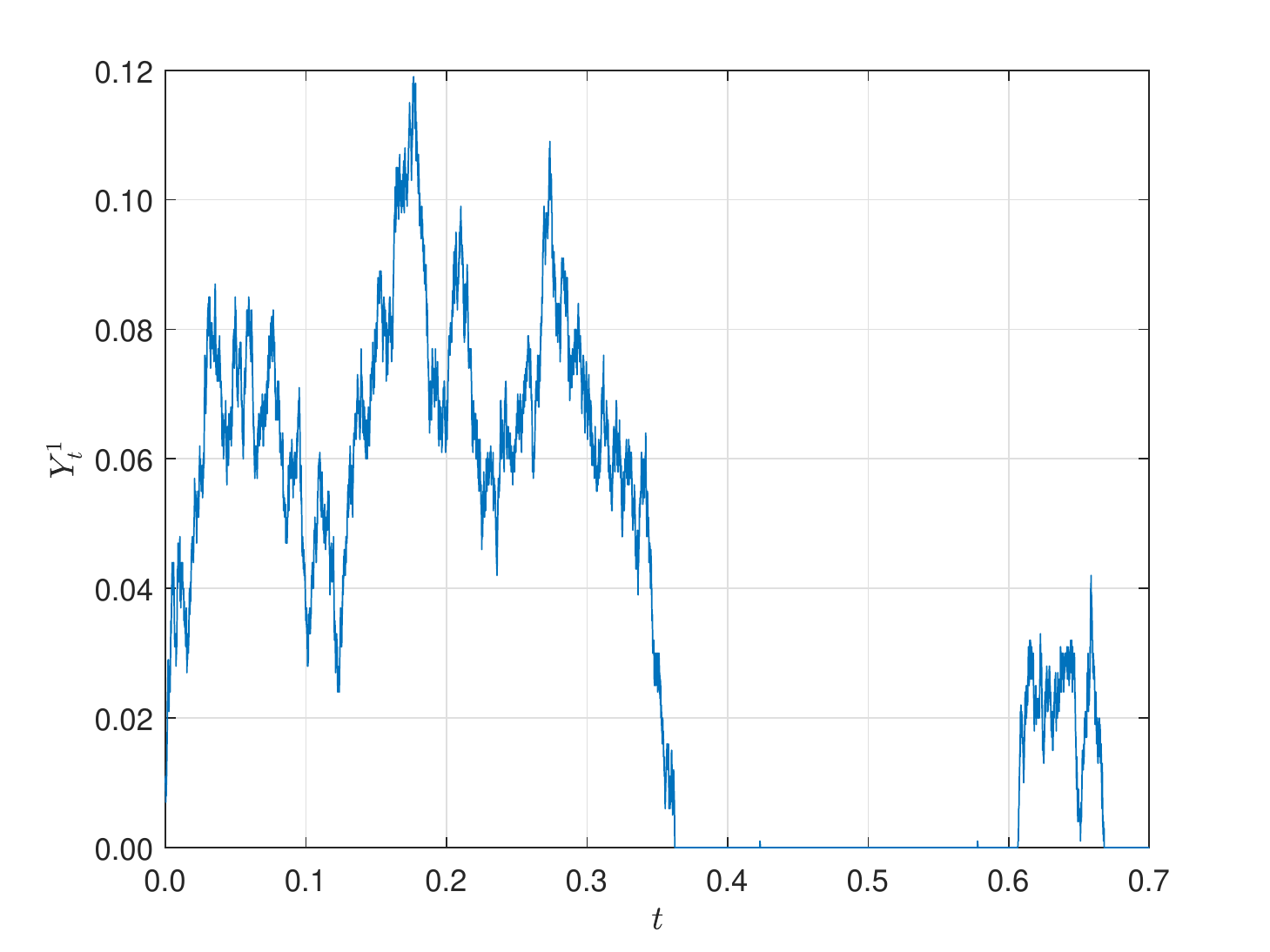}
	\end{subfigure}
	\begin{subfigure}{.496\textwidth}
		\centering
		\includegraphics[width=\linewidth]{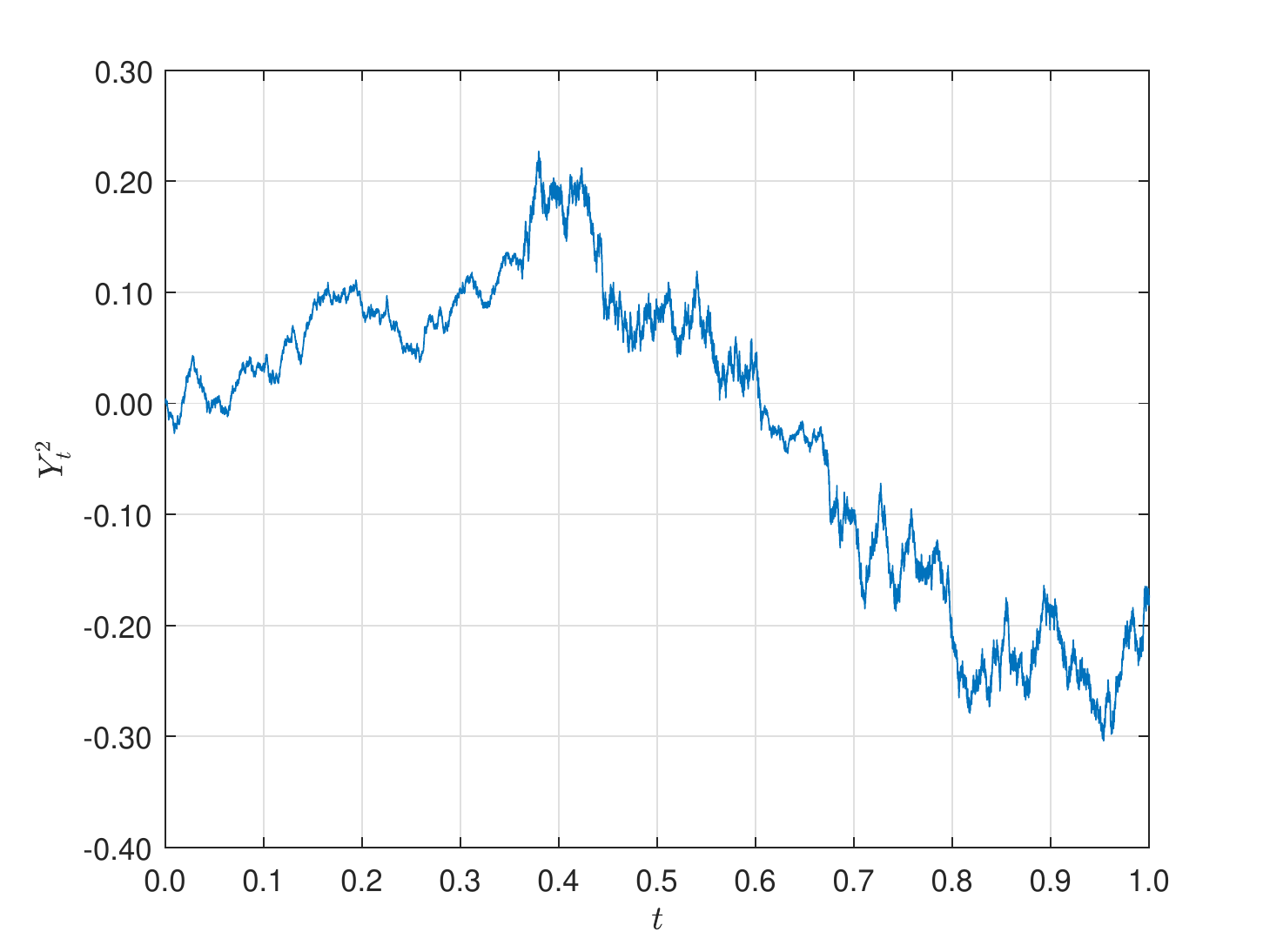}
	\end{subfigure}
	\begin{subfigure}{.496\textwidth}
		\centering
		\includegraphics[width=\linewidth]{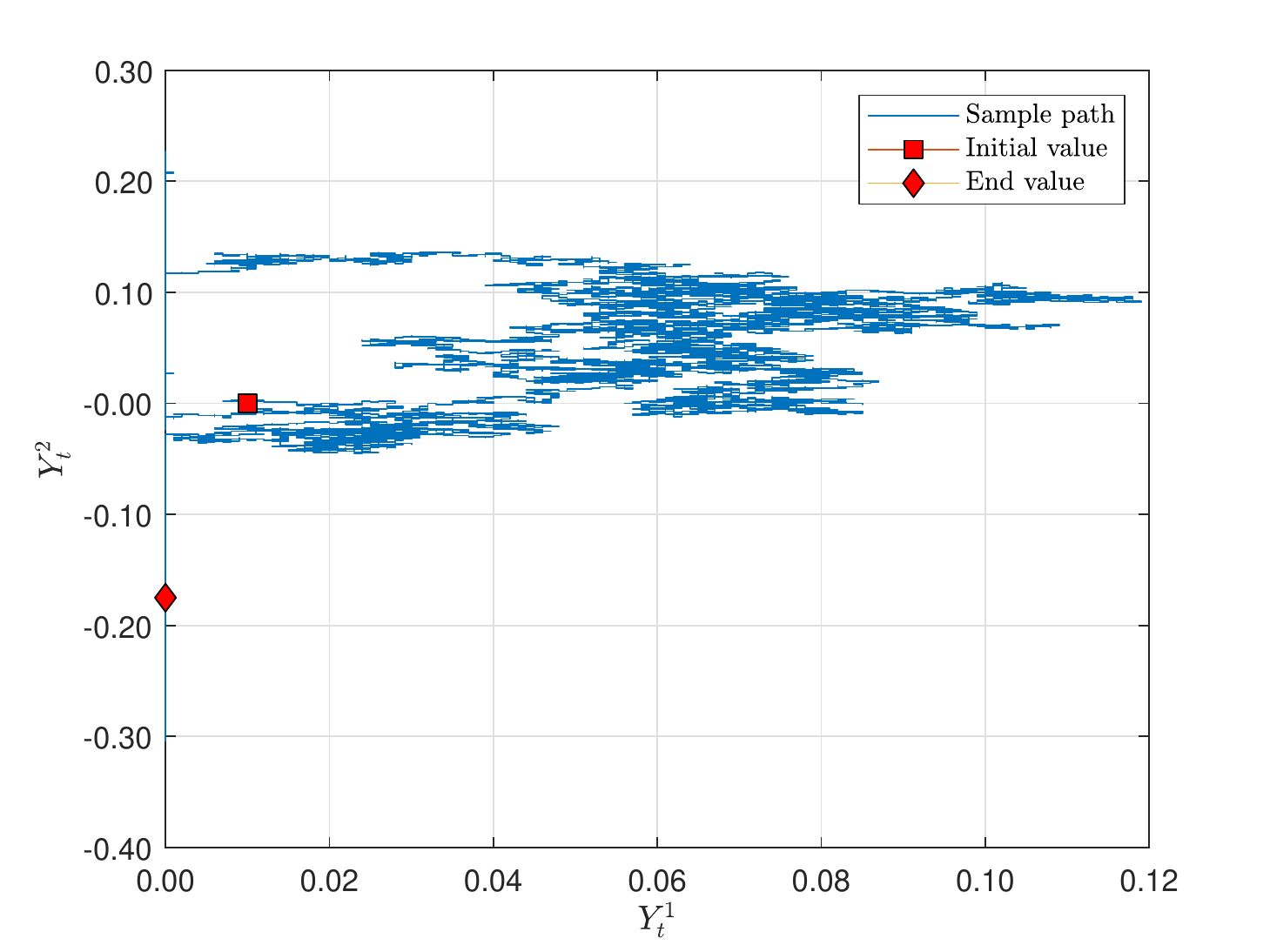}
	\end{subfigure}
	\caption{One sample path over the time interval $[0,1]$ generated from the CTMC with $h=1/1000$ constructed by the eigendecomposition approach for the sticky short rate model \eqref{eq:2D_interest_rate_model}. The first row shows the two dimensions separately while the second row shows the movement of the process in the two dimensional state space.}
	\label{fig:2D_sticky_model}
\end{figure}

Figure \ref{fig:2D_sticky_model} displays a sample path of the CTMC approximating the sticky OU process.
One can see that the first dimension exhibits stickiness at zero whereas the second dimension is unbounded. Both dimensions exhibit mean reverting behavior as implied by the dynamics given in \eqref{eq:2D_interest_rate_model}.
Furthermore, the starting point is in the interior and so, one can see that the process initially behaves like a two-dimensional OU process.
Once the first coordinate hits zero, the stickiness at the boundary is determined by the value of the second coordinate. The function $\nu(X^2)$ is increasing in $X^2$ and it shows the drift of $X^1$ to leave the boundary while it is there. Thus, the smaller $X^2$, the larger the stickiness for $X^1$ at the boundary. For the CTMC, the first coordinate sticks at zero for quite a long time as the second coordinate is quite negative, although short excursions into the interior of the state space might happen for the original diffusion.

\subsubsection{Convergence Results}
\label{subsec:2D_interest_rate_model_results}
\begin{figure}[htbp!]
	\centering
	\begin{subfigure}{.496\textwidth}
		\centering
				\includegraphics[width=\linewidth]{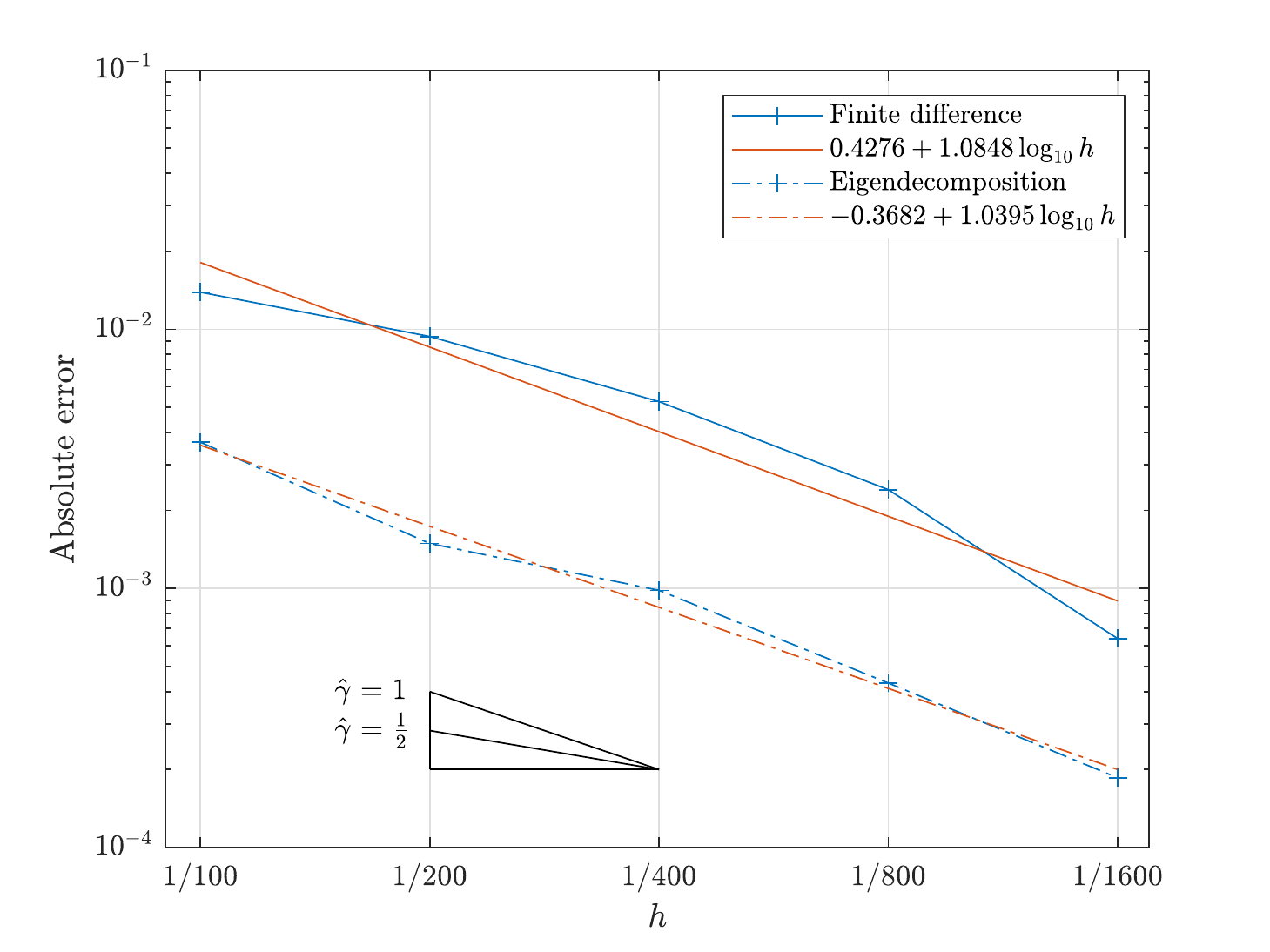}
	\end{subfigure}
	\begin{subfigure}{.496\textwidth}
		\centering
				\includegraphics[width=\linewidth]{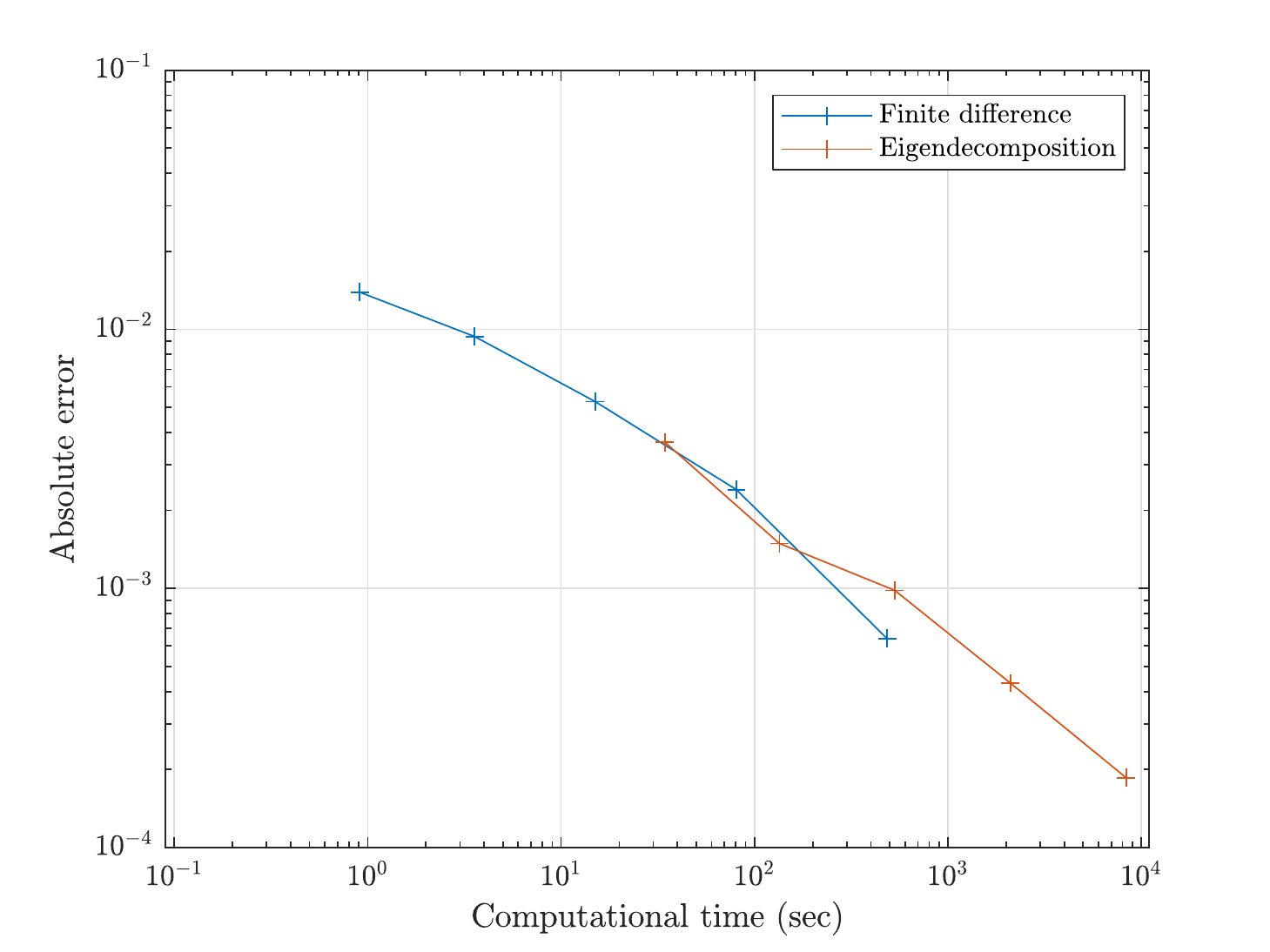}
	\end{subfigure}
	\caption{Convergence rate (left) and absolute error vs. computational time (right) for CTMC simulation of the sticky short rate model \eqref{eq:2D_interest_rate_model}. Both plots are on log-log scale. Each marker in both plots corresponds to one level of $h$.}
	\label{fig:convergence_2D_ir_model}
\end{figure}

To show convergence of our method for this model, we again set $h$ according to \eqref{eq:h}. For each value, we generate $10^5$ paths from the CTMC to approximate $v(T,x)$. Figure \ref{fig:convergence_2D_ir_model} shows the results for the finite difference and eigendecomposition approach using exact simulation for the CTMC. From the left plot, the theoretical convergence order is again verified. Moreover, for the same level of $h$, the eigendecomposition approach is more accurate.

The right plot displays absolute error with computational time. For error levels between $10^{-2}$ and $10^{-3}$, these two methods take roughly the same amount of time. For a fixed level of $h$, the eigendecomposition approach is slower due to additional computations related to adjusting the step size.

\section{Conclusion}
\label{sec:conclusion}

We develop a new simulation method for multidimensional diffusions with sticky boundaries, in which we simulate from a CTMC that approximates the diffusion. We provide two approaches for constructing the CTMC: the finite difference and the eigendecomposition approach. The second approach has two theoretical advantages. First, it always guarantees valid transition rates while the first approach cannot for problems with strong correlations. Second, for the same step size $h$, the error of the second approach is noticeably smaller thanks to the better transition directions it provides. 

Although directly computing the transition probability of a high-dimensional CTMC suffers from the curse of dimensionality due to the exponentially growing number of states, simulation of the process is free from the curse. As our construction shows, the total rate moving out of a state is $O(d/h^2)$. Consequently, exact simulation of the CTMC involves $O(d/h^2)$ number of transitions over a unit time interval. Moreover, the cost of simulating one transition is $O(d^2)$ for the finite difference approach or $O(d)$ for the eigendecomposition approach. As for additional calculations, their costs are also polynomials of $d$. Therefore, the total cost is a polynomial of $d$, which is not exponentially growing. If we simulate the CTMC by discretizing time, the number of transitions over a unit time interval is only $O(1/h)$, which does not depend on $d$ and the computational cost is further reduced.

We demonstrate the performance of our method in two nontrivial examples, from which we show that Monte Carlo simulation based on our method can provide accurate results for the value function and the generated sample paths can capture the sticky behavior. 

\section*{Acknowledgements}
The research of Lingfei Li was supported by Hong Kong Research Grant Council General Research Fund Grant 14202117. 
The research of Gongqiu Zhang was supported by National Natural Science Foundation of China Grant 11801423 and Shenzhen Basic Research Program Project JCYJ20190813165407555

\appendix

\section{Proofs}
\label{sec:proofs}

\noindent\emph{Proof of Theorem \ref{th:generator_multidimensional_sticky_diffusion}}: Using Lemma 2.11 in \S5 of Chapter 2 in \cite{dynkin1965}, the strong and weak infinitesimal generator of $\mathcal{P}_t$ coincide, i.e., we can calculate the limit in \eqref{eq:generator-def} pointwisely. For $f\in C_0^2(\bar{\mathbb{S}})$, Ito's formula shows that
\begin{align*}
	df\left(X_t\right)&=(\partial_xf\left(X_t\right))^\top\left(I\left(X_t\in\mathbb{S}\right)\mu\left(X_t\right)+I\left(X_t\in\partial\mathbb{S}\right)\hat{\beta}\left(X_t\right)\right)dt \nonumber\\
	&+(\partial_xf\left(X_t\right))^\top\left(I\left(X_t\in\mathbb{S}\right)\Sigma\left(X_t\right)dB_{1,t}+I\left(X_t\in\partial\mathbb{S}\right)\hat{\Gamma}\left(X_t\right)dB_{2,t}\right) \nonumber\\
	&+\frac{1}{2}\left(I\left(X_t\in\mathbb{S}\right)\textrm{Tr}\left(\Sigma\left(x\right)^\top\left(\partial_{xx}f\right)\Sigma\left(x\right)\right)+I\left(X_t\in\partial\mathbb{S}\right)\textrm{Tr}\left(\hat{\Gamma}\left(x\right)^\top\left(\partial_{xx}f\right)\hat{\Gamma}\left(x\right)\right)\right)dt.
\end{align*}
Applying this result yields
\begin{align*}
	\mathcal{G}f\left(x\right)&=\underset{t\searrow 0}{\lim}\ \frac{\mathbb{E}_x\left(f\left(X_t\right)\right)-f\left(x\right)}{t}= \underset{t\searrow 0}{\lim}\ \frac{\mathbb{E}_x\left(\int_0^td\left(f\left(X_s\right)\right)\right)}{t} \\
	&=I\left(X_t\in\mathbb{S}\right)\mathcal{A}f(x) + I\left(X_t\in\partial\mathbb{S}\right)\mathcal{K}f(x).
\end{align*}
The Wentzell boundary condition on $f$ results from the requirement that $\mathcal{G}f\in C_0(\bar{\mathbb{S}})$.\qed

\bigskip
\noindent\emph{Proof of Proposition \ref{prop:proper_transitions_diagonally_dominant}}: We first note that
\begin{equation}\label{eq:a-rates}
	2a_{i,\pm he_i}(x)h^2=A^{i,i}\left(x\right)-\sum_{j=1,j\neq i}^d\left\vert A^{i,j}\left(x\right)\right\vert-\mu^i\left(x\right)h.
\end{equation}
Set 
\begin{equation}
	\bar{h}=\underset{i=1,\ldots,d,\ x\in\mathbb{S}}{\min}\ \frac{A^{i,i}\left(x\right)-\sum_{i,j=1,j\neq i}^d\left\vert A^{i,j}\left(x\right)\right\vert}{\left\vert\mu^i\left(x\right)\right\vert},
\end{equation}
which is positive from the strict diagonal dominance assumption. It follows that for $h\leq \bar{h}$, 
\begin{equation}\label{eq:transition_rates_eq4}
	A^{i,i}\left(x\right)-\sum_{j=1,j\neq i}^d\left\vert A^{i,j}\left(x\right)\right\vert-\left\vert\mu^i\left(x\right)\right\vert h\geq 0,
\end{equation}
for all $i=1,\ldots,d$. Consequently, we obtain $a_{i,\pm he_i}\geq 0$ from \eqref{eq:a-rates}.
\qed

\bigskip
\noindent\emph{Proof of Theorem \ref{th:convergence_rate_ctmc}}: We first note that Theorem 3.1 in \cite{zeng1994} implies that $v\in C^{2,1}(\mathcal{S})$.

Denote the state space of the CTMC by $\bar{\mathbb{S}}^h$.
Let $Q(\cdot,\cdot)$ be the transition kernel of the CTMC defined by its transition rates $a$. Specifically, $Q(x,\{y\})$ is the rate of transitioning from $x$ to $y$ with $Q(x,\{x\})=0$ and $Q(x,E)=\sum_{y\in E}Q(x,\{y\})$ where $E\subseteq\bar{\mathbb{S}}^h$ (see Section 7.2 in \cite{durrett1996}). Let $\mathcal{P}^h_t$ be the transition operator of the CTMC $Y_t$, i.e., 
$\mathcal{P}^h_tf\left(x\right)=\mathbb{E}_x\left(f\left(Y_t\right)\right)$. The transition semigroup $(\mathcal{P}^h_t)_{t\ge0}$ is a strongly continuous semigroup of contractions on $C_0(\bar{\mathbb{S}}^h)$. Theorem 2.1 in Chapter 7 of \cite{durrett1996} shows that the infinitesimal generator of the CTMC is given by
\begin{align}
	\mathcal{Q}f\left(x\right)&=\int_{\bar{\mathbb{S}}^h}Q\left(x,dy\right)\left(f\left(y\right)-f\left(x\right)\right) \nonumber\\
	&=\sum_{y\in\bar{\mathbb{S}}^h}Q\left(x,\{y\}\right)f\left(y\right)-f\left(x\right)Q\left(x,\bar{\mathbb{S}}^h\right),\label{eq:generator_ctmc}
\end{align}
with $\mathcal{D}(\mathcal{Q})=\mathcal{D}(\mathcal{P}^h)=C_0(\bar{\mathbb{S}}^h)$.

Consider the error $e(T,x)$ defined by
\begin{equation*}
	e\left(T,x\right)=\mathbb{E}_x\left(f\left(Y_T\right)\right)-\mathbb{E}_x\left(f\left(X_T\right)\right),\qquad x\in\bar{\mathbb{S}}^h.
\end{equation*}
Using Lemma 6.2 in Chapter 1 of \cite{ethier2005}, we obtain for $x\in\bar{\mathbb{S}}^h$, 
\begin{align*}
	e\left(T,x\right)&=\mathbb{E}_x\left(f\left(Y_T\right)\right)-\mathbb{E}_x\left(f\left(X_T\right)\right)=\mathcal{P}^h_Tf\left(x\right)-\mathcal{P}_Tf\left(x\right) \\
	&=\mathcal{P}^h_Tf\left(x\right)-f\left(x\right)-\mathcal{P}_Tf\left(x\right)+f\left(x\right) \\
	&=\int_0^T\mathcal{P}^h_{T-t}\left(\mathcal{Q}-\mathcal{G}\right)\mathcal{P}_tf\left(x\right)dt.
\end{align*}
Consequently, we have
\begin{align*}
	\left\Vert e\left(T,x\right)\right\Vert_{\infty}&\leq\int_0^T\left\Vert\mathcal{P}^h_{T-t}\left(\mathcal{Q}-\mathcal{G}\right)\mathcal{P}_tf\left(x\right)\right\Vert_{\infty} dt \\
	&\leq\int_0^T\left\Vert\left(\mathcal{Q}-\mathcal{G}\right)\mathcal{P}_tf\left(x\right)\right\Vert_{\infty} dt,
\end{align*}
as $\mathcal{P}^h_{T-t}$ is a contraction and $\|\cdot\|_{\infty}$ is the maximum norm over $\bar{\mathbb{S}}^h$.

We next estimate the difference $\left(\mathcal{Q}-\mathcal{G}\right)v(t,x)$ for $x\in\bar{\mathbb{S}}^h$ for the two approaches separately. We use $\delta(x)$ as the adjusted step size for point $x$ (see \eqref{eq:definition_delta_x}). Although $\delta(x)$ depends on the transition direction in general, we do not reflect this dependence in the notation for simplicity. 

(I) Finite difference approach: For given $t\in(0,T]$ and $x\in\bar{\mathbb{S}}^h$, there holds now
\begin{align*}
	&\left(\mathcal{Q}-\mathcal{G}\right)v\left(t,x\right) \\
	&=\sum_{i=1}^dI\left(x\in\mathbb{S}_{\textrm{loc}}\right)\mu^i\left(x\right)\frac{v\left(t,x+\delta(x)e_i\right)-v\left(t,x-\delta(x)e_i\right)}{2\delta(x)}-I\left(x\in\mathbb{S}_{\textrm{loc}}\right)\mu^{i}\left(x\right)\frac{\partial}{\partial x^i}v\left(t,x\right)\\
	&+\frac{1}{2}\sum_{i=1}^dI\left(x\in\mathbb{S}_{\textrm{loc}}\right)A^{i,i}\left(x\right)\frac{v\left(t,x+\delta(x)e_i\right)-2v\left(t,x\right)+v\left(t,x-\delta(x)e_i\right)}{\delta\left(x\right)^2} \\
	&-\frac{1}{2}\sum_{i=1}^dI\left(x\in\mathbb{S}_{\textrm{loc}}\right)A^{i,i}\left(x\right)\frac{\partial^2}{\partial x^i\partial x^i}v\left(t,x\right) \\
	&+\sum_{i,j=1,j\neq i}^dI\left(x\in\mathbb{S}_{\textrm{loc}}\right)\bigg[I\left(A^{i,j}\left(x\right)\geq 0\right)A^{i,j}\left(x\right)\bigg(\frac{2v\left(t,x\right)-v\left(t,x+\delta(x)e_i\right)}{2\delta(x)^2} \\
	&+\frac{-v\left(t,x-\delta(x)e_i\right)-v\left(t,x+\delta(x)e_j\right)-v\left(t,x-\delta(x)e_j\right)}{2\delta(x)^2} \\
	&+\frac{v\left(t,x+\delta(x)e_i+\delta(x)e_j\right)+v\left(t,x-\delta(x)e_i-\delta(x)e_j\right)}{2\delta(x)^2}\bigg)\\
	&+I\left(A^{i,j}\left(x\right)<0\right)A^{i,j}\left(x\right)\bigg(\frac{-2v\left(t,x\right)+v\left(t,x+\delta(x)e_i\right)+v\left(t,x-\delta(x)e_i\right)}{2\delta(x)^2} \\
	&+\frac{v\left(t,x+\delta(x)e_j\right)+v\left(t,x-\delta(x)e_j\right)}{2\delta(x)^2} \\
	&+\frac{-v\left(t,x-\delta(x)e_i+\delta(x)e_j\right)-v\left(t,x+\delta(x)e_i-\delta(x)e_j\right)}{2\delta(x)^2}\bigg)\bigg] \\
	&-\frac{1}{2}\sum_{i,j=1,j\neq i}^dI\left(x\in\mathbb{S}_{\textrm{loc}}\right)A^{i,j}\left(x\right)\frac{\partial^2}{\partial x^i\partial x^j}v\left(t,x\right) \\
	&+\sum_{i\in\mathcal{I}_0(x)}I\left(x\in\partial\mathbb{S}_{\textrm{loc}}\right)\hat{\beta}^i\left(x\right)\frac{v\left(t,x+\delta(x)e_i\right)-v\left(t,x\right)}{\delta(x)}\\
	&+\sum_{i\in\mathcal{I}_0(x)^c}I\left(x\in\partial\mathbb{S}_{\textrm{loc}}\right)\hat{\beta}^i\left(x\right)\frac{v\left(t,x+\delta(x)e_i\right)-v\left(t,x-\delta(x)e_i\right)}{2\delta(x)}\\
	&-\sum_{i=1}^dI\left(x\in\partial\mathbb{S}_{\textrm{loc}}\right)\hat{\beta}^i\left(x\right)\frac{\partial}{\partial x^i}v\left(t,x\right)\\
	&+\frac{1}{2}\sum_{i\in\mathcal{I}_0(x)^c}I\left(x\in\partial\mathbb{S}_{\textrm{loc}}\right)\hat{G}^{i,i}\left(x\right)\frac{v\left(t,x+\delta(x)e_i\right)-2v\left(t,x\right)+v\left(t,x-\delta(x)e_i\right)}{\delta\left(x\right)^2} \\
	&-\frac{1}{2}\sum_{i\in\mathcal{I}_0(x)^c}I\left(x\in\partial\mathbb{S}_{\textrm{loc}}\right)\hat{G}^{i,i}\left(x\right)\frac{\partial^2}{\partial x^i\partial x^i}v\left(t,x\right) \\
	&+\sum_{i,j\in\mathcal{I}_0(x)^c,j\neq i}I\left(x\in\partial\mathbb{S}_{\textrm{loc}}\right)\bigg[I\left(\hat{G}^{i,j}\left(x\right)\geq 0\right)\hat{G}^{i,j}\left(x\right)\bigg(\frac{2v\left(t,x\right)-v\left(t,x+\delta(x)e_i\right)}{2\delta(x)^2} \\
	&+\frac{-v\left(t,x-\delta(x)e_i\right)-v\left(t,x+\delta(x)e_j\right)-v\left(t,x-\delta(x)e_j\right)}{2\delta(x)^2} \\
	&+\frac{v\left(t,x+\delta(x)e_i+\delta(x)e_j\right)+v\left(t,x-\delta(x)e_i-\delta(x)e_j\right)}{2\delta(x)^2}\bigg)\\
	&+I\left(\hat{G}^{i,j}\left(x\right)<0\right)\hat{G}^{i,j}\left(x\right)\bigg(\frac{-2v\left(t,x\right)+v\left(t,x+\delta(x)e_i\right)+v\left(t,x-\delta(x)e_i\right)}{2\delta(x)^2} \\
	&+\frac{v\left(t,x+\delta(x)e_j\right)+v\left(t,x-\delta(x)e_j\right)}{2\delta(x)^2} \\
	&+\frac{-v\left(t,x-\delta(x)e_i+\delta(x)e_j\right)-v\left(t,x+\delta(x)e_i-\delta(x)e_j\right)}{2\delta(x)^2}\bigg)\bigg] \\
	&-\frac{1}{2}\sum_{i,j\in\mathcal{I}_0(x)^c,j\neq i}I\left(x\in\partial\mathbb{S}_{\textrm{loc}}\right)\hat{G}^{i,j}\left(x\right)\frac{\partial^2}{\partial x^i\partial x^j}v\left(t,x\right).
\end{align*}

Application of Taylor's theorem yields now
\begin{align*}
	&\frac{v\left(t,x+\delta\left(x\right)e_i\right)-v\left(t,x-\delta\left(x\right)e_i\right)}{2\delta\left(x\right)}\\
	&=\frac{\partial}{\partial x^i}v\left(t,x\right)+\frac{1}{2}\frac{\partial^2}{\partial x^i\partial x^i}v\left(t,x+\xi_1\delta\left(x\right)e_i\right)\delta\left(x\right) -\frac{1}{2}\frac{\partial^2}{\partial x^i\partial x^i}v\left(t,x-\xi_2\delta\left(x\right)e_i\right)\delta\left(x\right)\ (\xi_1,\xi_2\in[0,1])\\
	&=\frac{\partial}{\partial x^i}v\left(t,x\right)+O(\delta^2(x)).
\end{align*}
For the last equality, we use the Lipschitz continuity of the $\partial_{xx}v$.

Similar derivations can be done for the second order derivative resulting in
\begin{align*}
	&\frac{v\left(t,x+\delta\left(x\right)e_i\right)-2v\left(t,x\right)+v\left(t,x-\delta\left(x\right)e_i\right)}{\delta\left(x\right)^2}=\frac{\partial^2}{\partial x^i\partial x^i}v\left(t,x\right)+O\left(\delta\left(x\right)\right).
\end{align*}
Likewise, the approximation error of the cross derivative terms is also of order $O(\delta(x))$ ($i,j=1,\ldots,d$ and $j\neq i$).
At the boundary, we obtain for $i\in\mathcal{I}_0(x)$:
\begin{align*}
	\frac{v\left(t,x+\delta\left(x\right)e_i\right)-v\left(t,x\right)}{\delta\left(x\right)}&=\frac{\partial}{\partial x^i}v\left(t,x\right)+O(\delta\left(x\right)) \\
	\frac{v\left(t,x\right)-v\left(t,x-\delta\left(x\right)e_i\right)}{\delta\left(x\right)}&=\frac{\partial}{\partial x^i}v\left(t,x\right)+O(\delta\left(x\right)).
\end{align*}
Aggregating these estimates, we obtain $\left(\mathcal{Q}-\mathcal{G}\right)v\left(t,x\right)=O(\delta(x))$.

(II) Eigendecomposition based approach: let $u_i(x),u_i^{\hat{G}}(x)\in\mathbb{R}^d$ with $\Vert u_i(x)\Vert=1$ be the $i$-th normalized eigenvector of $A$ for $i=1,\ldots,d$ and $u_i^{\hat{G}}(x)$ the normalized eigenvector of $\hat{G}$. 
The difference between the two infinitesimal generators can now be written in the following way: for $x\in\bar{\mathbb{S}}^h$,
\begin{align*}
	&\left(\mathcal{Q}-\mathcal{G}\right)v\left(t,x\right) \\
	&=I\left(x\in\mathbb{S}_{\textrm{loc}}\right)\frac{v\left(t,x+\delta(x)\mu\left(x\right)\right)-v\left(t,x\right)}{\delta(x)}-\sum_{i=1}^dI\left(x\in\mathbb{S}_{\textrm{loc}}\right)\mu^i\left(x\right)\frac{\partial}{\partial x^i}v\left(t,x\right)\\
	&+\frac{1}{2}\sum_{i=1}^dI\left(x\in\mathbb{S}_{\textrm{loc}}\right)\lambda_i\left(x\right)\frac{v\left(t,x+\delta(x)u_i\left(x\right)\right)-2v\left(t,x\right)+v\left(t,x-\delta(x)u_i\left(x\right)\right)}{\delta\left(x\right)^2} \\
	&-\frac{1}{2}\sum_{i,j=1}^dI\left(x\in\mathbb{S}_{\textrm{loc}}\right)A^{i,j}\left(x\right)\frac{\partial^2}{\partial x^i\partial x^j}v\left(t,x\right) \\
	&+I\left(x\in\partial\mathbb{S}_{\textrm{loc}}\right)\frac{v\left(t,x+\delta(x)\hat{\beta}\left(x\right)\right)-v\left(t,x\right)}{\delta(x)}-\sum_{i=1}^dI\left(x\in\partial\mathbb{S}_{\textrm{loc}}\right)\hat{\beta}^i\left(x\right)\frac{\partial}{\partial x^i}v\left(t,x\right)\\
	&+\frac{1}{2}\sum_{i=1}^dI\left(x\in\partial\mathbb{S}_{\textrm{loc}}\right)\lambda_i^{\hat{G}}\left(x\right)\frac{v\left(t,x+\delta(x)u_i^{\hat{G}}\left(x\right)\right)-2v\left(t,x\right)+v\left(t,x-\delta(x)u_i^{\hat{G}}\left(x\right)\right)}{\delta\left(x\right)^2} \\
	&-\frac{1}{2}\sum_{i,j=1}^dI\left(x\in\partial\mathbb{S}_{\textrm{loc}}\right)\hat{G}^{i,j}\left(x\right)\frac{\partial^2}{\partial x^i\partial x^j}v\left(t,x\right).
\end{align*}
Note that from the eigendecomposition, there holds $\sum_{i=1}^d\lambda_i\left(x\right)u_i^n\left(x\right)\left(u_i^n\left(x\right)\right)^\top=A\left(x\right)$ and so
\begin{align*}
	&\frac{1}{2}\sum_{i=1}^d\lambda_i\frac{v\left(t,x+\delta\left(x\right)u_i\left(x\right)\right)-2v\left(t,x\right)+v\left(t,x-\delta\left(x\right)u_i\left(x\right)\right)}{\delta\left(x\right)^2} \\
	&=\frac{1}{2}\sum_{i=1}^d\lambda_i\left(x\right)u_i\left(x\right)^\top\partial_{xx}v\left(t,x\right)u_i\left(x\right)+O\left(\delta\left(x\right)\right) \\
	&=\frac{1}{2}\sum_{i=1}^d\lambda_i\left(x\right)\sum_{k=1}^d\sum_{l=1}^du_i^{k}\left(x\right)\left(\partial_{xx}v\left(t,x\right)\right)^{k,l}u_i^{l}\left(x\right)+O\left(\delta\left(x\right)\right) \\
	&=\frac{1}{2}\sum_{k=1}^d\sum_{l=1}^dA^{k,l}\left(x\right)\frac{\partial^2}{\partial x^k\partial x^l}v\left(t,x\right)+O\left(\delta\left(x\right)\right).
\end{align*}
The terms involving $\mu$ and $\hat{\beta}$ can easily be handled and by their boundedness, the error is of first order. 
Moreover, using the eigendecomposition of $\hat{G}$, the approximation error of the second order derivative is also $O(\delta(x))$.

(III) Using the results in (I) and (II) together with $\delta(x)\leq h$, we obtain
\begin{align*}
	\left\Vert e\left(T,x\right)\right\Vert_{\infty}&\leq\int_0^T\left\Vert\left(\mathcal{Q}-\mathcal{G}\right)\mathcal{P}_tf\left(x\right)\right\Vert_{\infty} dt \\
	&\leq\int_0^T O(h) dt \\
	&\leq CTh
\end{align*}
for some constant $C>0$. This shows that the CTMC approximation converges with first order. \qed
\bibliographystyle{apalike}
\bibliography{references}

\end{document}